\documentclass[opre,nonblindrev]{informs3} 

\OneAndAHalfSpacedXI 

\usepackage{endnotes}
\let\footnote=\endnote
\let\enotesize=\normalsize
\def\notesname{Endnotes}%
\def\makeenmark{$^{\theenmark}$}
\def\enoteformat{\rightskip0pt\leftskip0pt\parindent=1.75em
  \leavevmode\llap{\theenmark.\enskip}}

\usepackage{algorithm2e}
\RestyleAlgo{ruled}
\usepackage{diagbox}
\usepackage{multirow}
\usepackage[symbol]{footmisc}

\usepackage[hidelinks]{hyperref}
\usepackage{xcolor}
\def\boxit#1{%
  \smash{\color{black}\fboxrule=1pt\relax\fboxsep=2pt\relax%
  \llap{\rlap{\fbox{\vphantom{0}\makebox[#1]{}}}~}}\ignorespaces
}
\usepackage{array}
\newcolumntype{C}[1]{>{\centering\let\newline\\\arraybackslash\hspace{0pt}}m{#1}}

\TheoremsNumberedThrough     
\ECRepeatTheorems

\EquationsNumberedThrough    

\begin{document}


\RUNAUTHOR{Li, Li, and Zhang}

\RUNTITLE{High-Probability Guarantees for Stochastic Saddle Point Problems}

\TITLE{General Procedure to Provide High-Probability Guarantees for Stochastic Saddle Point Problems}

\ARTICLEAUTHORS{%
\AUTHOR{Dongyang Li, Haobin Li, Junyu Zhang\footnote{Corresponding author.}}
\AFF{Department of Industrial Systems Engineering and Management, National University of Singapore, Singapore\\ \{\EMAIL{dongyang\_li@u.nus.edu}, \EMAIL{li\_haobin@nus.edu.sg}, \EMAIL{junyuz@nus.edu.sg}\}}
} 

\ABSTRACT{%
This paper considers smooth strongly convex and strongly concave (SC-SC) stochastic saddle point (SSP) problems. Suppose there is an arbitrary oracle that \emph{in expectation} returns an $\epsilon$-solution in the sense of certain gaps, which can be the duality gap or its weaker variants. We propose a general PB-SSP framework to guarantee an $\epsilon$ small duality gap solution \emph{with high probability} via only $\mathcal{O}\big(\log \frac{1}{p}\cdot\text{poly}(\log \kappa)\big)$ calls of this oracle, where $p\in(0,1)$ is the confidence level and $\kappa$ is the condition number. When applied to the sample average approximation (SAA) oracle, in addition to equipping the solution with high probability, our approach even improves the sample complexity by a factor of $\text{poly}(\kappa)$, since the high-probability argument enables us to circumvent some key difficulties of the uniform stability analysis of SAA.
}%



\maketitle

%

\section{Introduction}
In this paper, we consider the stochastic saddle point problem
\begin{equation}
\label{equ:minmax}
    \min_{x \in \mathcal{X}} \max_{y \in \mathcal{Y}} \Phi(x,y):=\mathbb{E}\left[\Phi_{\xi}(x,y) \right],
\end{equation}
where $\mathcal{X}$ and $\mathcal{Y}$ are closed and convex sets, and $\xi$ is a random variable satisfying an unknown distribution $\mathcal{P}$. This formulation finds a wide range of applications
in adversarial learning \cite{goodfellow2020generative,sinha2017certifying}, reinforcement learning \cite{puterman2014markov,wang2017primal}, robust optimization \cite{namkoong2017variance}, and game theory \cite{von2007theory,roughgarden2010algorithmic}, and so on. We will focus on the basic setting where $\Phi_\xi(\cdot,\cdot)$ is smooth, convex in $x$ and concave in $y$ for almost every $\xi\sim\mathcal{P}$. For any feasible solution $(\hat x,\hat y)\in\mathcal{X}\times\mathcal{Y}$, we denote the duality gap as 
\begin{equation*}
    \Delta_{\Phi}(\hat{x},\hat{y}):=\max_{y\in\mathcal{Y}} \Phi(\hat{x},y)-\min_{x\in\mathcal{X}} \Phi(x,\hat{y}).
\end{equation*}
The standard results on convex-concave (C-C) SSP problems usually upper bound the total number of required samples for finding solutions with $\epsilon$ small \emph{expected} duality gap, recent examples include stochastic approximation (SA) type algorithms \cite{nemirovski2009robust,shalev2013stochastic,zhang2017stochastic,yan2019stochastic} and the sample average approximation (SAA) method \cite{zhang2021generalization}, to name a few. Despite the bulk of literature that guarantees the small duality gap in expectation, the high-probability results are rare. 
Now, let us consider a weaker variant of the duality gap: 
\begin{equation*}
    \Delta^w_{\Phi}(\hat{x},\hat{y}):=\Phi(\hat{x},y^*)-\Phi(x^*,\hat{y}),
\end{equation*}
where $(x^*,y^*)$ is the optimal solution to (\ref{equ:minmax}). For $\mu_x$-strongly convex and $\mu_y$-strongly concave saddle point problem, it is not hard to see that 
$$\frac{\mu_x}{2}\|\hat{x}-x^*\|^2 + \frac{\mu_y}{2}\|\hat{y}-y^*\|^2\leq \Delta^w_{\Phi}(\hat{x},\hat{y})\leq \Delta_{\Phi}(\hat{x},\hat{y}).$$
Therefore, for SC-SC saddle point problems, given an arbitrary oracle that returns a solution $(\hat x,\hat y)\in\mathcal{X}\times\mathcal{Y}$ such that   
$$\mathbb{E}[\Delta_\Phi(\hat x,\hat y)]\leq \epsilon\qquad\mbox{or}\qquad \mathbb{E}[\Delta_\Phi^w(\hat x,\hat y)]\leq \epsilon.$$
The goal of this research is to produce a solution $(\bar x,\bar y)\in\mathcal{X}\times\mathcal{Y}$ such that 
\begin{equation}
    \mathbb{P}\left[\Delta_{\Phi}(\bar x,\bar y) \leq \epsilon \right] \geq 1-p 
    \label{equ: high prob guarantee}
\end{equation}
via only a few executions of this oracle, where $p\in(0,1)$ is the confidence level. 

Indeed, there are few naive approaches to achieve \eqref{equ: high prob guarantee} for general oracles, without relying on the additional light-tail noise assumptions. The first method is to apply the oracle to generate a solution $(\hat x,\hat y)$ such that $\mathbb{E}[\Delta_\Phi(\hat x,\hat y)]\leq p\cdot\epsilon$, then \eqref{equ: high prob guarantee} is immediately guaranteed by Markov's inequality. However, the sample complexity to compute such a $p\epsilon$-accurate solution is often proportional to $(p\epsilon)^{-1}$ for SC-SC problems. Naively applying Markov's inequality will amplify the sample complexity by an addition $\mathcal{O}\big(\frac{1}{p}\big)$ factor. The second approach is to generate a group of solutions $\{(\hat x_i,\hat y_i)\}_{i=1}^m$ with $\epsilon$ expected duality gap. Then if one can evaluate the $\Delta_\Phi(\hat x_i,\hat y_i)$ to $\epsilon$ accuracy, then the pair $(\hat x_{i^*},\hat y_{i^*})$ with minimal estimated duality gap will satisfy \eqref{equ: high prob guarantee} as long as $m=\Omega\big(\ln \frac{1}{p}\big)$. This resolves the unfavorable dependence on $\mathcal{O}\big(\frac{1}{p}\big)$, yet evaluating the duality gap for each solution will itself consume $\mathcal{O}(\epsilon^{-2})$ samples, which is much more expensive than the usual $\tilde{\mathcal{O}}(\epsilon^{-1})$ sample complexity of the oracle itself. It is also possible to apply the robust distance estimation technique for unconstrained SC-SC SSP problems, but that will often bring in additional $\mathcal{\kappa}$ factors in the sample complexity, see detailed discussion in Section \ref{sec: preliminary}. 

Therefore, we would like to propose a general meta-framework to generate solutions satisfying \eqref{equ: high prob guarantee} while only suffering an $\mathcal{O}\big(\ln \frac{1}{p}\cdot\text{poly}(\ln \kappa)\big)$ overhead in the sample complexity. The proposed framework generalizes the previous \underline{P}rox\underline{B}oost algorithm \cite{davis2021low} for minimization problems to the more general \underline{S}tochastic \underline{S}addle \underline{P}oint problems, thus we call the framework PB-SSP.

\subsection{Main contribution}\label{subsection: main contribution}
Given any SSP oracle that bounds the duality gap $\Delta_{\Phi}(\hat{x},\hat{y})$ or its weaker variant $\Delta^w_{\Phi}(\hat{x},\hat{y})$ in expectation, we provide a general meta framework called PB-SSP to equip this oracle with high-probability guarantees in the sense of (\ref{equ: high prob guarantee}), while only increasing sample complexity by a factor of $\mathcal{O}\big(\log \frac{1}{p}\cdot\text{poly}(\log \kappa)\big)$. More formally, for an arbitrary SSP oracle $\mathcal{M}$, we denote $\mathcal{M}^w(\Phi,\delta)$ as the operation of calling $\mathcal{M}$ to generate a solution $(\hat{x},\hat{y})$ s.t. $\mathbb{E}\left[\Delta^w_{\Phi}(\hat{x},\hat{y})\right] \leq \delta$, and we denote the corresponding sample complexity as $C_{\mathcal{M}}^w(\Phi,\delta)$. Similarly, we denote $\mathcal{M}(\Phi,\delta)$ as the operation of generating $(\hat{x},\hat{y})$ s.t. $\mathbb{E}\left[\Delta_{\Phi}(\hat{x},\hat{y})\right] \leq \delta$, and the corresponding sample complexity as $C_{\mathcal{M}}(\Phi,\delta)$. Apparently, for the same oracle $\mathcal{M}$, we have $C_{\mathcal{M}}^w(\Phi,\delta) \leq C_{\mathcal{M}}(\Phi,\delta)$. In particular, when $\mathcal{M}$ is the SAA oracle \cite{zhang2021generalization}, we have $C_{\mathcal{M}}(\Phi,\delta) = \kappa \cdot C^w_{\mathcal{M}}(\Phi,\delta)$. 
By operating $\mathcal{M}^w(\cdot)$ for a sequence of perturbed proximal point subproblems, the PB-SSP framework proposed in this paper will return 
a point $(\bar{x},\bar{y})$ satisfying (\ref{equ: high prob guarantee}) with a sample complexity of
\begin{equation*}
    \ln\left(\frac{\ln(\kappa)}{p}\right)\ln(\kappa)\cdot C_{\mathcal{M}}^w\left(\Phi,\frac{\epsilon}{\ln(\kappa)}\right).
\end{equation*}

An interesting application of our PB-SSP framework is the sample average approximation (SAA) oracle, which is frequently used in practice. In the SAA oracle, the user is provided with some offline dataset $\Gamma:=\{\xi_1,\xi_2,\cdots,\xi_n\}$ drawn i.i.d. from $\mathcal{P}$. Then SAA generates the output by solving
\begin{equation}
\label{equ: saa minmax}
    (\hat{x},\hat{y})=\argmin_{x \in \mathcal{X}}\,\, \argmax_{y \in \mathcal{Y}} \hat{\Phi}_n(x,y)\,:=\,\frac{1}{n}\sum_{i=1}^n \Phi_{\xi_i}(x,y),
\end{equation} 
regardless of the detailed saddle point algorithms. Existing work \cite{zhang2021generalization} analyzes the SAA oracle by leveraging the uniform stability arguments. For SSP problems with SC-SC modulus $\mu>0$ and function Lipschitz constant $\ell$, an $\mathcal{O}\big(\frac{\ell^2}{n\mu}\big)$ bound on $\mathbb{E}\big[\Delta^w_{\Phi}(\hat{x},\hat{y})\big]$ is proved. With addition $L$-Lipschitz assumption on the gradients, an $\mathcal{O}\big(\kappa\cdot\frac{\ell^2}{n\mu}\big)$ bound on $\mathbb{E}\big[\Delta_{\Phi}(\hat{x},\hat{y})\big]$ can be guaranteed. By incorporating the PB-SSP framework, we are able to circumvent some key difficulties in the uniform stability analysis and yield an $\tilde{\mathcal{O}}\big(\frac{\ell^2}{n\mu}\big)$ bound on $\Delta_{\Phi}(\hat{x},\hat{y})$ with high probability, where the condition number $\kappa$ only appears in the logarithmic factors. Compared to the original SAA oracle, the PB-SSP not only provides a high-probability guarantee but also improves the bound by a factor of $\kappa$. This improvement can be huge if $\kappa\gg 1$. In particular, when the SC-SC modulus comes from the  $+\frac{\mu}{2}\|x\|^2-\frac{\mu}{2}\|y\|^2$ regularization added to some general convex-concave problem, the PB-SSP can improve the bound on the duality gap from $\mathcal{O}(1/\sqrt[3]{n})$ in expectation to $\tilde{\mathcal{O}}(1/\sqrt{n})$ with high probability. We summarize the improvement for SAA oracle in Table \ref{table: complexity bounds}.\newpage
\begin{table}[htb]
\centering
\caption{Sample complexity bounds for SAA oracles on SSP problems.}
\begin{tabular}{l | c c | c}
\hline
     & \multicolumn{2}{c|}{SAA \cite{zhang2021generalization}}  & SAA + PB-SSP (this paper)\\
\hline
    \multirow{2}{*}{\backslashbox{Problems}{Bounds}} & \multicolumn{2}{c|}{\quad Expectation guarantees $\quad$ } & high-probability guarantees\\
    & $\Delta_{\Phi}^w$ & $\Delta_{\Phi}$ & $\Delta_{\Phi}$\\
\hline
    SC-SC unconstrained & $\mathcal{O}\left(\frac{C\kappa^4}{\mu\epsilon}\right)$ & $\mathcal{O}\left(\frac{C\kappa^4}{\mu\epsilon}\right)$ & $\mathcal{O}\left(\ln(\kappa)\ln\left(\frac{\ln(\kappa)}{p}\right)
    \frac{C\kappa^2}{\mu \epsilon}
    \right)$\\
    SC-SC constrained & $\mathcal{O}\left(\frac{\ell^2}{\mu\epsilon}\right)$ & $\mathcal{O}\left(\frac{\ell^2\kappa}{\mu\epsilon}\right)$ & $\mathcal{O}\left(\ln^2(\kappa)\ln\left(\frac{\ln(\kappa)}{p}\right)
    \frac{\ell^2}{\mu \epsilon}
    \right)$\\
    C-C constrained & $\mathcal{O}\left(\frac{\ell^2 D^2}{\epsilon^2}\right)$ & $\mathcal{O}\left(\frac{\ell^2L D^4}{\epsilon^3}\right)^\dagger$ & $\mathcal{O}\left(\ln^2\left(\frac{LD^2}{\epsilon}\right)\ln\left(\frac{\ln(LD^2/\epsilon)}{p}\right)
    \frac{\ell^2D^2}{\epsilon^2}
    \right)$\\
\hline 
\end{tabular}  
\label{table: complexity bounds}
\end{table}
\vspace{-0.9cm}
 {\footnotesize $\quad\,\dagger\,\,$This bound is not provided in \cite{zhang2021generalization}, we derived it based on the analysis of \cite[Theorem 2 \& Theorem 3]{zhang2021generalization}}.

\subsection{Related works}
Currently, most existing SSP algorithms belong to stochastic approximation (SA) approach. Let $n$ be the number of samples consumed by the algorithm. When the SSP problem is only convex-concave, to our best knowledge, Nemirovski and Rubinstein \cite{nemirovski2002efficient} are the first to establish an $\mathcal{O}(1/\sqrt{n})$ bound on duality gap by extending Polyak’s method \cite{polyak1990new}. Similar bounds are also obtained for the stochastic mirror descent ascent algorithm in \cite{nemirovski2009robust}. Chen et al. \cite{chen2014optimal} considered a class of SSP problems with a bilinear coupling term and leveraged Nesterov's smoothing \cite{nesterov2005smooth} scheme to accelerate the primal-dual methods. Their work retains the $\mathcal{O}\left(1/\epsilon^2\right)$ sample complexity for ensuring an $\epsilon$ small duality gap while improving the dependence on several parameters. Zhao \cite{zhao2022accelerated} considered SSP problems with a three-composite structure and also reached the $\mathcal{O}(1/\sqrt{n})$ bound. Besides, Beznosikov et al. \cite{beznosikov2020gradient} developed zeroth-order methods for general non-smooth C-C SSP problems. 
Under a more general setting, SSP problems fall into a special case of stochastic variational inequality (SVI) problems. When the operator is monotone, $\mathcal{O}(1/\sqrt{n})$ bound can be established. For example, Juditsky et al. \cite{juditsky2011solving} proposed the stochastic mirror-prox algorithm, and Chen et al. \cite{chen2017accelerated} accelerated it for a special class of SVI problems. More recent works on monotone SVI problems include \cite{iusem2017extragradient,gidel2018variational,bach2019universal,hsieh2020explore,mishchenko2020revisiting,gorbunov2022clipped}. 

When it comes to SC-SC SSP (or strong monotone SVI) problems, tighter bounds can be expected. For example, given $n$ samples, \cite{natole2018stochastic,hsieh2019convergence,fallah2020optimal,huang2022new,gorbunov2022clipped} derived $\mathcal{O}(1/n)$ bounds for the squared distance from the saddle point (i.e., $\|\hat{x}-x^*\|^2+\|\hat{y}-y^*\|^2$). Yan et al. \cite{yan2019stochastic} and Yan et al. \cite{yan2020optimal} provided $\mathcal{O}(1/n)$ bound for the primal gap (i.e., $\max_{y\in\mathcal{Y}} \Phi(\hat{x},y)-\Phi(x^*,y^*)$) and duality gap respectively. There also exist some works considering saddle point problems with finite-sum structure and applying stochastic primal-dual methods, e.g., see \cite{shalev2013stochastic,palaniappan2016stochastic,zhang2017stochastic,chavdarova2019reducing,luo2019stochastic,alacaoglu2022stochastic}. It is worth noting the finite-sum structures are often the direct result of sample average approximation of the stochastic objective functions. Therefore, these methods can be considered as the subroutines for solving the deterministic empirical problem constructed by SAA. 

Finally, concerning the SAA approach for SSP or SVI problems, the research is extremely limited. Xu \cite{xu2010sample} investigated the SAA solution of general SVI problems and showed the asymptotic convergence without any finite sample analysis. Recently, Zhang et al. \cite{zhang2021generalization} utilized the uniform stability-based approach and bounded the $\mathbb{E}\big[\Delta^w_{\Phi}(\hat{x},\hat{y})\big]$ by $\mathcal{O}(1/\sqrt{n})$ for C-C SSP problems and the $\mathbb{E}\big[\Delta_{\Phi}(\hat{x},\hat{y})\big]$ by $\mathcal{O}(1/n)$ for SC-SC problems respectively. Lei et al. \cite{lei2021stability} established stability bounds for gradient-based algorithms when applied to the SAA counterpart of SSP problems. Similar stability results of gradient-based algorithms appear in \cite{farnia2021train}.

There are also very few high-probability results for SSP problems in the literature, and most of them heavily rely on the sub-Gaussian assumption for the stochastic gradient. For example, to obtain a solution with $\epsilon$ small duality gap with high probability, Nemirovski et al. \cite{nemirovski2009robust} bounded the sample complexity by
$\mathcal{O}\left(\ln^2\left(\frac{1}{p}\right)\frac{1}{\epsilon^2}\right)$ for C-C SSP problems. Juditsky et al. \cite{juditsky2011solving} derived similar results for monotone SVI problems. Yan et al. \cite{yan2020optimal} proposed an Epoch-GDA algorithm for SC-SC cases and derived $\mathcal{O}\left(\ln\left(\frac{\ln(1/\epsilon)}{p}\right)\frac{1}{ \epsilon}\right)$ sample complexity bound to guarantee the duality gap is $\epsilon$ small. Laguel et al. \cite{laguel2023high} investigate the accelerated primal-dual algorithm for SC-SC problems and ensure $\epsilon$ small squared distance with sample complexity of $\mathcal{O}\left(\ln\left(\frac{1}{p}\right)\ln\left(\frac{1}{\epsilon}\right)\frac{1}{\epsilon}\right)$. Chen et al. \cite{chen2014optimal,chen2017accelerated} and Zhao \cite{zhao2022accelerated} also gave high-probability results for various SSP problems with special structures like bilinear coupling and smooth components. Under stronger assumption that the stochastic gradients are almost surely bounded, Yan et al. \cite{yan2019stochastic} investigate a special class of SC-SC SSP problems where the primal-dual coupling term is linear in the dual variable. It proposes a restart scheme and provides high-probability guarantees for the primal gap.  All the results mentioned above need the sub-Gaussian or even stronger assumptions, except for Gorbunov et al. \cite{gorbunov2022clipped} that derived the high-probability results for unconstrained SVI problems, where the authors modified several first-order methods with the gradient clipping technique. As a result, high-probability guarantees are established for the duality gap and the squared distance in the monotone and strong monotone cases respectively.
Compared with \cite{gorbunov2022clipped} which is unconstrained and algorithm-dependent, our procedure establishes a more general framework and allows us to handle constrained problems.

\section{Preliminaries} \label{sec: preliminary}
\textbf{Basic Assumptions.} $ $ First, let us formally state a few fundamental assumptions of the objective function. Unless otherwise stated, they are always assumed throughout this paper.
\begin{assumption}
\label{assump:SC-SC}
$\exists$ $\mu_x, \mu_y >0$ s.t. for almost every $\xi \sim \mathcal{P}$ and $\forall x\in\mathcal{X}, \forall y\in\mathcal{Y}$,  the function $\Phi_{\xi}(\cdot,y)$ is $\mu_x$-strongly convex and the function $\Phi_{\xi}(x,\cdot)$ is $\mu_y$-strongly concave. Namely,
\begin{equation*}
\begin{aligned}
    &\Phi_{\xi}(x_2,y)\geq \Phi_{\xi}(x_1,y)+\langle \nabla_x\Phi_{\xi}(x_1,y),x_2-x_1 \rangle +\frac{\mu_x}{2}\|x_1-x_2\|^2, \quad \forall x_1,x_2 \in \mathcal{X}, y\in \mathcal{Y},\\
    &\Phi_{\xi}(x,y_2)\leq \Phi_{\xi}(x,y_1)+\langle \nabla_y\Phi_{\xi}(x,y_1),y_2-y_1 \rangle -\frac{\mu_y}{2}\|y_1-y_2\|^2,\,\,\, \quad \forall y_1,y_2 \in \mathcal{Y}, x\in \mathcal{X}.\\
\end{aligned}
\end{equation*}
Consequently, we say $\Phi$ is SC-SC with modulus $(\mu_x, \mu_y)$.
\end{assumption}

\begin{assumption}
\label{assump:smoothness}
$\exists$ $L_x, L_y, L_{xy}>0$ s.t. for 
$\forall x_1,x_2 \in \mathcal{X}$ and $\forall y_1,y_2 \in \mathcal{Y}$, we have
\begin{equation*}
\begin{aligned}
    &\|\nabla_x \Phi(x_1,y_1)-\nabla_x\Phi(x_2,y_1)\|\leq L_x\|x_1-x_2\|, \quad
    \|\nabla_y \Phi(x_1,y_1)-\nabla_y\Phi(x_1,y_2)\|\leq L_y\|y_1-y_2\|,\\
    &\|\nabla_x \Phi(x_1,y_1)-\nabla_x\Phi(x_1,y_2)\|\leq L_{xy}\|y_1-y_2\|, \quad
    \|\nabla_y \Phi(x_1,y_1)-\nabla_y\Phi(x_2,y_1)\|\leq L_{xy}\|x_1-x_2\|.
\end{aligned}
\end{equation*}
\end{assumption}

\begin{assumption}
\label{assump: function Lipschitz}
The feasible regions $\mathcal{X}$ and $\mathcal{Y}$ are closed and convex. $\exists$ $\ell_x,\ell_y>0$ s.t. for almost every $\xi\sim\mathcal{P}$, it holds that 
\begin{equation*}
\begin{aligned}
    |\Phi_{\xi}(x_2,y)-\Phi_{\xi}(x_1,y)| \leq \ell_x\|x_1-x_2\|\quad\mbox{and}\quad 
    |\Phi_{\xi}(x,y_1)-\Phi_{\xi}(x,y_2)| \leq \ell_y\|y_1-y_2\|.
\end{aligned}
\end{equation*}
for $\forall x_1,x_2,x \in \mathcal{X} \text{ and } \forall y_1,y_2,y \in \mathcal{Y}$.
\end{assumption}
In the case of SC-SC SSP problems, Assumption \ref{assump: function Lipschitz} typically applies when the feasible regions are compact. In the upcoming Section \ref{sec: unconstrained}, we will first address the simpler setting with unbounded domains, and Assumption \ref{assump: function Lipschitz} will be temporarily dropped in this setting.
The above and subsequent $\|\cdot\|$ denotes the Euclidean norm. To simplify notation, we further denote $\mu:=\min\{\mu_x,\mu_y\}$, $L:=\max\{L_x,L_y,L_{xy}\}$, $\ell:=\max\{\ell_x,\ell_y\}$. The condition number of $\Phi$ is denoted by $\kappa:=L/\mu$.\vspace{0.5cm}

\noindent\textbf{Robust Distance Estimation with Pseudometrics.} $ $
The robust distance estimation (RDE) \cite{nemirovskij1983problem,hsu2016loss} is the key tool of the PB-SSP framework. Consider any pseudometric $\rho: \mathcal{X}\times\mathcal{X} \rightarrow \mathbb{R}$ that is non-negative, symmetric, and satisfies the triangle inequality, that is, $\rho(x_1,x_2)\geq0$, $\rho(x_1,x_1)=0$, $\rho(x_1,x_2)=\rho(x_2,x_1)$, and $\rho(x_1,x_2)\leq \rho(x_1,x_3)+\rho(x_3,x_2)$, for $\forall x_1,x_2,x_3\in\mathcal{X}$. Denote $B_r^{\rho}(x)=\{y \in \mathcal{X}: \rho(x,y) \leq r\}$ the $r$-radius ball around $x$ under pseudometric $\rho$. Then the robust distance estimation technique extracts the candidate ``centers'' from a group of points by Algorithm \ref{alg:extract}.
\begin{algorithm}[htb!]
\SetKwInOut{Return}{Return}
\caption{Extract($\{x_j\}_{j=1}^m,\rho$)}\label{alg:extract}
\KwIn{A list of $m$ points $X=\{x_1,\dots,x_m\} \subset \mathcal{X}$, a pseudometric $\rho$ on $\mathcal{X}$.}
\textbf{for} $j=0,\dots,m$ \textbf{compute}: $r_j=\min\{r\geq0:|B_r^{\rho}(x_j) \cap X|>\frac{m}{2}\}$. \\
Compute the median radius $\hat{r}=\text{median}(r_1,\dots,r_m)$.\\
\Return{$\mathcal{I}=\{k \in [1,m]:r_k \leq \hat{r}\}$.}
\end{algorithm}
 With slight abuse, $\text{median}(r_1,...,r_m)$ in Algorithm \ref{alg:extract} denotes the $\left \lceil \frac{m}{2} \right \rceil$'th entry in the ordered list, so the algorithm will return at least $m/2$ points. Lemma \ref{lemma:robust distance estimation} illustrates how the RDE boosts the confidence for a collection of candidate points within a pseudometric space.
\begin{lemma}\textbf{\emph{(\cite[Lemma 11]{davis2021low})}}
\label{lemma:robust distance estimation}
Let $\rho$ and $X=\{x_1,\dots,x_m\}$ be the input of Algorithm \ref{alg:extract}. Suppose the point $x^*$ satisfies $\mathbb{P}[\rho(x_j,x^*)\leq \delta] \geq 2/3$ for some $\delta>0$, $\forall j$.
Then w.p. at least $1-\exp\left(-\frac{m}{18}\right)$, it holds that the event $E:=\{|B_{\delta}^{\rho}(x^*) \cap X| > \frac{m}{2}\}$ happens, and $E$ ensures
$\rho(x_k,x^*)\leq 3\delta$ for all $k\in\textnormal{Extract}(X,\rho)$.
\end{lemma}\vspace{0.5cm}

\noindent\textbf{The Primal and Dual Functions.} $ $ Back to the saddle point problem \eqref{equ:minmax}, for the objective function $\Phi(x,y)$, we define its primal function as $f(x):=\max_{y\in\mathcal{Y}}\Phi(x,y)$ and the dual function as $g(y):=\min_{x\in \mathcal{X}}\Phi(x,y)$. For SC-SC saddle point problems, the strong duality indicates that $f(x^*)=g(y^*)=\Phi(x^*,y^*)$. It is commonly known that (e.g., see Lemma A.5 of \cite{nouiehed2019solving}),
$f(x)$ is $\mu_x$-strongly convex and $L_f$-smooth, and $g(y)$ is $\mu_y$-strongly concave and $L_g$-smooth, where
\begin{equation}
\label{equ: L of f and g}
    L_{f} := L_x+L^2_{xy}/\mu_y \quad \text{and} \quad
    L_{g} := L_y+L^2_{xy}/\mu_x.
\end{equation}
Then, we have the following two-sided bounds for any $(x,y) \in \mathcal{X}\times \mathcal{Y}$. 
\begin{eqnarray}
\label{equ: constrained_two_side-1}
    \Delta^w_{\Phi}(x,y) & = & \Phi(x,y^*)-\Phi(x^*,y^*)+\Phi(x^*,y^*)-\Phi(x^*,y)\\
    & \overset{(a)}{\geq} & \frac{\mu_x}{2}\|x-x^*\|^2 + \langle\nabla_x \Phi(x^*,y^*) ,x-x^*\rangle
    +\frac{\mu_y}{2}\|y-y^*\|^2-\langle\nabla_y \Phi(x^*,y^*) ,y-y^*\rangle,\nonumber
\end{eqnarray}
and 
\begin{eqnarray}
\label{equ: constrained_two_side-2}
    \Delta_{\Phi}(x,y) & = & f(x)-f(x^*)+g(y^*)-g(y)\\
    &\overset{(b)}{\leq} & \frac{L_f}{2} \|x-x^*\|^2
    +\langle\nabla_x \Phi(x^*,y^*) ,x-x^*\rangle+\frac{L_g}{2}\|y-y^*\|^2
    -\langle\nabla_y \Phi(x^*,y^*) ,y-y^*\rangle.\nonumber
\end{eqnarray}
The step (a) is due to the SC-SC property in Assumption \ref{assump:SC-SC}, and the step (b) is due to the Lipschitz smoothness of $\nabla f(x)$ and $\nabla g(y)$ in (\ref{equ: L of f and g}) and the Danskin's theorem which implies $\nabla f(x^*)=\nabla_x \Phi(x^*,y^*)$ and $\nabla g(y^*)=\nabla_y \Phi(x^*,y^*)$. In particular, for unconstrained problems where $\mathcal{X}$ and $\mathcal{Y}$ are the full Euclidean spaces, both $\nabla_x \Phi(x^*,y^*)$ and $\nabla_y \Phi(x^*,y^*)$ are equal to 0. 

Note that for unconstrained saddle point problems where $\nabla_x \Phi(x^*,y^*)=\nabla_y \Phi(x^*,y^*)=0$, the duality gap is only upper bounded by the squared distance. This implies that, with the above inequalities and the robust distance estimation technique, one is already able to provide a high-probability solution without the extra overhead of $\mathcal{O}(1/p)$ or $\mathcal{O}(1/\epsilon)$ for unconstrained problems. Firstly, call $m=\lceil18\ln(1/p)\rceil$ times of $\mathcal{M}^w(\Phi,\delta)$ to obtain a sequence of candidate points $\{\hat{x}_j,\hat{y}_j\}_{j=1}^m$. Then (\ref{equ: constrained_two_side-1}) together with Markov's inequality implies $\mathbb P[\|\hat{x}_j-x^*\|\leq \sqrt{6\delta/\mu}]\geq 2/3$ and $\mathbb P[\|\hat{y}_j-y^*\|\leq \sqrt{6\delta/\mu}]\geq 2/3$, $\forall j$. Secondly, set $\rho$ to be Euclidean norm and use Algorithm $\ref{alg:extract}$ to get $\mathcal{I}_x=\text{Extract}(\{\hat{x}_j\}_{j=1}^m,\rho)$ and $\mathcal{I}_y=\text{Extract}(\{\hat{y}_j\}_{j=1}^m,\rho)$. Next, pick arbitrary $k_1\in \mathcal{I}_x$ and $k_2\in \mathcal{I}_y$ and set $(\bar{x},\bar{y})=(\hat{x}_{k_1},\hat{y}_{k_2})$. If the problem is unconstrained, by Lemma \ref{lemma:robust distance estimation} and applying the upper bound in (\ref{equ: constrained_two_side-2}), we can ensure
$\mathbb{P}\left[\Delta_{\Phi}(\bar{x},\bar{y})\leq 54(\kappa^2+\kappa)\delta \right]\geq 1-p$, then setting $\delta=\epsilon/54(\kappa^2+\kappa)$ provides the final result. The total computation cost will be 
\begin{equation}
    \mathcal{O}\left(\ln\left(\frac{1}{p}\right)\cdot C_{\mathcal{M}}^w\left(\Phi,\frac{\epsilon}{\kappa^2}\right)\right).
    \label{equ: robust distance bound}
\end{equation}
We name this plausible way to generate solutions with high-probability guarantees as the RDE approach in the paper. Since $C_{\mathcal{M}}^w\left(\Phi,\delta\right)$ is often of order $\tilde{\mathcal{O}}(1/\delta)$ for SC-SC problems, the RDE approach introduces an extra $\kappa^2$ factor to the overall computation cost. Moreover, the RDE approach is not directly applicable to the constrained problem due to the gradient terms in \eqref{equ: constrained_two_side-2}. This key drawback motivates us to adopt the inexact proximal point method in \cite{davis2021low}.\vspace{0.5cm} 

\noindent\textbf{Inexact Proximal Point Algorithm.} $ $ 
The inexact proximal point algorithm (IPPA) is a classical optimization method. For any objective function $h(\cdot)$, given an initial proximal center $x_0^c$, IPPA approximately solves a sequence of proximal point subproblems:
\begin{equation*}
    \min h^i(x):=h(x)+\frac{\lambda_i}{2}\|x-x_i^c\|^2, \quad \mbox{for}\quad i=0,1,\cdots,T,
\end{equation*}
where the next proximal center $x_{i+1}^c$ is often chosen as some approximate minimizer of $h^i(x)$. By standard arguments (e.g., see Theorem 2 of \cite{davis2021low}), IPPA satisfies the following lemma.
\begin{lemma}
\label{lemma: inexact proximal point}
For arbitrary sequences of $\{\lambda_i\}_{i=0}^T$ and $\{x_i^c\}_{i=0}^T$, denote $x_{i+1}^*:=\argmin_x h^i(x)$ as the exact minimizer of $h^i(x)$ and default $x_0^*=x^*:=\argmin_x h(x)$. Then we have
\begin{equation*}
    h(x_{t+1}^c)-h(x^*)\leq h^t(x_{t+1}^c)-h^t(x^*_{t+1})+\sum_{i=0}^t \frac{\lambda_i}{2}\|x_i^c-x_i^*\|^2,\quad \mbox{for}\quad t=0,1,\dots,T.
\end{equation*}
\end{lemma}

\section{Unconstrained SSP problems}\label{sec: unconstrained}
A key idea of our procedure is to bound $\Delta_{\Phi}(\bar{x},\bar{y})$ by controlling $f(\bar{x})-f(x^*)$ and $g(y^*)-g(\bar{y})$ separately. Noting $f(x)$ is convex, and $g(y)$ is concave, the inexact proximal point method can be adapted to mitigate the dependence on $\kappa$.
Naturally, we extend the IPPA to SSP problems, which forms a pillar for our procedure. Fix an increasing sequence of penalties $\lambda_x^0,\dots,\lambda_x^T$ and $\lambda_y^0,\dots,\lambda_y^T$, and a sequence of centers $x_0^c,\dots,x_T^c$ and $y_0^c,\dots,y_T^c$. For $i=0,\dots,T$, we define a series of perturbed functions together with their solutions as
\begin{equation}
\begin{aligned}
    &\Phi_x^i(x,y):=\Phi(x,y)+\frac{\lambda_x^i}{2}\|x-x_i^c\|^2, \quad
    f^i(x):=\max_{y \in \mathcal{Y}} \Phi_x^i(x,y) =f(x)+\frac{\lambda_x^i}{2}\|x-x_i^c\|^2,\\
    &\Phi_y^i(x,y):=\Phi(x,y)-\frac{\lambda_y^i}{2}\|y-y_i^c\|^2, \quad
    g^i(y):=\min_{x \in \mathcal{X}} \Phi_y^i(x,y) =g(y)-\frac{\lambda_y^i}{2}\|y-y_i^c\|^2;
\end{aligned}
\label{equ: perturbed functions}\nonumber
\end{equation}
\begin{equation}
    (x^*_{i+1},y^*_{x,i+1}) := \argmin_{x\in\mathcal{X}}\,\,\argmax_{y\in\mathcal{Y}} \Phi_x^i(x,y), \quad (x^*_{y,i+1},y^*_{i+1}) := \argmax_{y\in\mathcal{Y}}\,\,\argmin_{x\in\mathcal{X}} \Phi_y^i(x,y).
    \label{equ: perturbed function solutions}
\end{equation}

We set $\lambda_x^{-1}=\lambda_y^{-1}=0$ throughout the paper such that $\Phi_x^{-1}=\Phi_y^{-1}=\Phi$, $x^*_{0}=x^*$, $y^*_{0}=y^*$. It is easy to verify that $\Phi_x^i(x,y)$ is SC-SC modulus $(\mu_x+\lambda_x^i,\mu_y)$, and it is Lipschitz smooth in terms of $(L_x+\lambda_x^i,L_y,L_{xy})$. The situation for $\Phi_y^i(x,y)$ is similar. Hence, we deduce that $f^i(x)$ is $(\mu_x+\lambda_x^i)$-strongly convex and $L_f^i$-smooth, and $g^i(y)$ is $(\mu_y+\lambda_y^i)$-strongly concave and $L_g^i$-smooth, where
\begin{equation}
    \label{equ: L of f^i and g^i}
    L_{f}^i := L_x+L^2_{xy}/\mu_y +\lambda_x^i \quad \text{and} \quad
    L_{g}^i := L_y+L^2_{xy}/\mu_x +\lambda_y^i.
\end{equation}
The two-sided bounds in (\ref{equ: constrained_two_side-1}) and (\ref{equ: constrained_two_side-2}) can also be easily extended to $\Phi_x^i$ and $\Phi_y^i$, $i=0,\dots,T$, by replacing the notation accordingly. Applying Lemma \ref{lemma: inexact proximal point} to $f$ and $g$ yields the following proposition.  
\begin{proposition}
Set $\lambda_x^{-1}=\lambda_y^{-1}=0$, $x_0^*=x^*:=\argmin_{x \in \mathcal{X}} f(x)$ and $y_0^*=y^*:=\argmax_{y \in \mathcal{Y}} f(y)$. $x_i^*$ and $y_i^*$ are defined by (\ref{equ: perturbed function solutions}). For $t=0,\dots,T$, the following estimate holds:
\begin{equation*}
\begin{aligned}
    f(x_{t+1}^c)-f(x^*) \leq f^t(x_{t+1}^c)-f^t(x^*_{t+1})+
    \sum_{i=0}^t \frac{\lambda_x^i}{2}\|x_i^c-x^*_i\|^2,\\
    g(y^*)-g(y_{t+1}^c) \leq g^t(y^*_{t+1})-g^t(y_{t+1}^c)+
    \sum_{i=0}^t \frac{\lambda_y^i}{2}\|y_i^c-y^*_i\|^2.
\end{aligned}
\end{equation*}
\label{prop: decomposition of f and g}
\end{proposition}
\vspace{-1.2cm}
Proposition \ref{prop: decomposition of f and g} provides guidance for bounding the duality gap using the inexact proximal point method. If the estimated solutions, $x_{T+1}^c$ and $y_{T+1}^c$, of the last pair of perturbed functions, $f^T(x)$ and $g^T(y)$, are returned as $(\bar{x},\bar{y})$, the overall duality gap, $\Delta_{\Phi}(\bar{x},\bar{y})=f(\bar{x})-g(\bar{y})$, can be upper bounded by the suboptimality in the last pair of perturbed functions and the errors incurred along the way over both streams.

\subsection{The PB-SSP framework} \label{subsec: PB-SSP}
We are now prepared to outline our method. For the ease of understanding, we shall start with the unconstrained SSP problems and then extend the methodology to constrained problems in Section \ref{sec: constrained}. In the unconstrained setting, both $\mathcal{X}$ and $\mathcal{Y}$ in problem (\ref{equ:minmax}) are full Euclidean spaces, and hence Assumption \ref{assump: function Lipschitz} will be dropped in this section due to the lack of domain compactness. The bounds in this section shall be derived under Assumptions \ref{assump:SC-SC} and \ref{assump:smoothness}. Note that in this setting, we have $\nabla_x \Phi(x^*,y^*)=\nabla_x \Phi(x^*,y^*)=0$, and the two-sided bounds (\ref{equ: constrained_two_side-1}) and (\ref{equ: constrained_two_side-2}) reduce to
\begin{equation}
\label{equ: unconstrained two-side}
    \frac{\mu_x}{2}\|x-x^*\|^2+ \frac{\mu_y}{2}\|y-y^*\|^2
    \leq \Delta^w_{\Phi}(x,y)
    \leq \Delta_{\Phi}(x,y)
    \leq \frac{L_f}{2}\|x-x^*\|^2+\frac{L_g}{2}\|y-y^*\|^2.
\end{equation}
To achieve a high-probability guarantee for the duality gap, we will monitor the overall function error and the total computation cost over a sequence of perturbed optimization problems. Since $\lambda_x^i$ and $\lambda_y^i$ are increasing, we may gradually decrease the tolerance on the errors $\|x_i^c-x_i^*\|$ and $\|y_i^c-y_i^*\|$, along with which the condition numbers $\frac{L_f+\lambda_x^i}{\mu_x+\lambda_x^i}$ and $\frac{L_g+\lambda_y^i}{\mu_y+\lambda_y^i}$ of the perturbed functions are decreasing to $\mathcal{O}(1)$. With this in mind, we introduce the PB-SSP algorithm as a general framework to boost the confidence of any SSP oracle $\mathcal{M}$.

Note that, this algorithm also depends on the undefined sequence $\{\lambda_x^i\}_{i=0}^T$ and $\{\lambda_y^i\}_{i=0}^T$. To simplify notation, we treat them as global parameters specified in theorems rather than as the algorithm input. The following Theorem \ref{theorem: pro boost of SSP} summarizes the guarantees of the PB-SSP algorithm.

\begin{algorithm}[H]
\SetKwInOut{Return}{Return}
\caption{PB-SSP($\delta,p,T$)}\label{alg:PB-SSP}
\KwIn{$\delta > 0$, $p \in (0,1)$, $T\in \mathbb{N}$.}
Set $\lambda_x^{-1}=\lambda_y^{-1}=0$, $x_{-1}^c=y_{-1}^c=0$.\\
\For{$i=0,\dots,T$}{
Set $\epsilon_x^{i}=\sqrt{\frac{2\delta}{\mu_x+\lambda_x^{i-1}}}$ and $\epsilon_y^{i}=\sqrt{\frac{2\delta}{\mu_y+\lambda_y^{i-1}}}$.\\
Generate a point $(x_{i}^c,y_{i}^c)$ satisfying
\begin{equation}
\begin{aligned}
    \mathbb{P}[\|x_{i}^c-x^*_{i}\|\leq\epsilon_x^{i}]\geq 1-\frac{p}{2T+4} \quad \text{and} \quad
    \mathbb{P}[\|y_{i}^c-y^*_{i}\|\leq\epsilon_y^{i}]\geq 1-\frac{p}{2T+4},
    \label{equ: high prob x and y}
\end{aligned} 
\end{equation}
where $x_i^* := \argmin_{x \in \mathcal{X}} f^{i-1}(x)$ and $y_i^* := \argmax_{y\in \mathcal{Y}}g^{i-1}(y)$.
}
Generate a point $(x_{T+1}^c,y_{T+1}^c)$ satisfying
\begin{equation}
    \mathbb{P}[f^T(x_{T+1}^c)-f^T(x^*_{T+1})\leq \delta ] \geq 1-\frac{p}{2T+4} \quad \text{and} \quad
    \mathbb{P}[g^T(y^*_{T+1})-g^T(y_{T+1}^c)\leq \delta ] \geq 1-\frac{p}{2T+4}.
    \label{equ: high prob x and y last round}
\end{equation}\\
\Return{$(x_{T+1}^c,y_{T+1}^c)$}
\end{algorithm}

\begin{theorem}[Proximal Boost of SSP]\label{theorem: pro boost of SSP}
Fix a target relative accuracy $\delta>0$, a probability of failure $p \in (0,1)$, and an integer $T>0$. Then with probability at least $1-p$, the point $(x_{T+1}^c,y_{T+1}^c)=\textnormal{PB-SSP}(\delta,p,T)$ satisfies
\begin{equation}
    \Delta_{\Phi}(x_{T+1}^c,y_{T+1}^c) \leq 
    \delta\left(2+\sum_{i=0}^T\frac{\lambda_x^i}{\mu_x+\lambda_x^{i-1}}+\frac{\lambda_y^i}{\mu_y+\lambda_y^{i-1}}\right).
    \label{equ: PB-SSP}
\end{equation}
\label{them: PB-SSP}
\end{theorem}
\proof{Proof}
Denote the events $E_x^i:=\{\|x_{i}^c-x^*_{i}\|\leq\epsilon_x^{i}\}$ and $E_y^i:=\{\|y_{i}^c-y^*_{i}\|\leq\epsilon_y^{i}\}$, $i=0,\dots,T$. Also denote $E_x^{T+1}:=\{f^T(x_{T+1}^c)-f^T(x^*_{T+1})\leq \delta\}$ and $E_y^{T+1}:=\{g^T(y^*_{T+1})-g^T(y_{T+1}^c)\leq \delta\}$. We deduce
\begin{equation*}
    \mathbb{P}\left[\bigcap_{i=0}^{T+1}\left(E_x^i\cap E_y^i\right)\right]=1-\mathbb{P}\left[\bigcup_{i=0}^{T+1}\left((E_x^i)^c\cup (E_y^i)^c\right)\right] \geq 1-\sum_{i=0}^{T+1}\left(\frac{p}{2T+4}+\frac{p}{2T+4}\right)=1-p.
\end{equation*}
Given the occurrence of the event $\left\{\bigcap_{i=0}^{T+1}\left(E_x^i\cap E_y^i\right)\right\}$, Proposition \ref{prop: decomposition of f and g} and (\ref{equ: unconstrained two-side}) indicates that 
\begin{equation*}
\begin{aligned}
    \Delta_{\Phi}(x_{T+1}^c,y_{T+1}^c)
    =&f(x_{T+1}^c)-f(x^*)+g(y^*)-g(y_{T+1}^c)\\
    \leq & f^T(x_{T+1}^c)-f^T(x^*_{T+1})+ g^T(x^*_{T+1})-g^T(y_{T+1}^c) +
    \sum_{i=0}^T \frac{\lambda_x^i}{2}\|x_i^c-x^*_i\|^2 + \frac{\lambda_y^i}{2}\|y_i^c-y^*_i\|^2\\
    \leq & \delta\left(2+\sum_{i=0}^T\frac{\lambda_x^i}{\mu_x+\lambda_x^{i-1}}+\frac{\lambda_y^i}{\mu_y+\lambda_y^{i-1}}\right),
\end{aligned}
\end{equation*}
where the last inequality uses the definition of $\epsilon_x^j$ and $\epsilon_y^j$. This completes the proof.
\Halmos
\endproof

Looking at (\ref{equ: PB-SSP}), we follow the choice from \cite{davis2021low} to let $\lambda_x^i=\mu_x \nu^i$ and $\lambda_y^i=\mu_y \nu^i$ for some constant base number $\nu$. After only $T=\mathcal{O}(\log_{\nu}\kappa)$ 
iterations, the condition numbers of $f^T$ and $g^T$, i.e., $\frac{L_f+\lambda_x^T}{\mu_x+\lambda_x^T}$ and $\frac{L_g+\lambda_y^T}{\mu_y+\lambda_y^T}$, will reduce to $\mathcal{O}(1)$. Therefore, based on Proposition \ref{prop: decomposition of f and g}, we only need to bound distances for the first $T+1$ iterations via standard robust distance estimation, while applying the RDE approach to bound the function gaps for the last iteration. But the difference is that, for the last iteration where the problem condition numbers are $\mathcal{O}(1)$, the dependence on squared condition numbers of RDE  (\ref{equ: robust distance bound}) becomes no longer an issue. Corollary \ref{corollary: PBSSP with Geo decay} provides parameter settings along with the sample complexity bound for achieving the desirable high-probability guarantee.

\begin{corollary}[Proximal Boost of SSP with geometric decay]\label{corollary: PBSSP with Geo decay}
Let us fix an arbitrary target accuracy $\epsilon>0$, and a probability of failure $p \in (0,1)$. Define the algorithm parameters:
\begin{equation*}
    \begin{aligned}
    &T=\left \lceil \log_{\nu}\left( \max\left(\frac{L_{xy}^2/\mu_y+L_x}{\mu_x},\frac{L_{xy}^2/\mu_x+L_y}{\mu_y}\right)
    \right) \right \rceil,\\
    &\delta=\frac{\epsilon}{4+4T}, \quad \lambda_x^i=\mu_x \nu^i \quad \text{and} \quad \lambda_y^i=\mu_y \nu^i \quad \text{with} \quad \nu=2 \quad \forall i \in \{0,1,\cdots,T\}.
\end{aligned}
\end{equation*}
Then the point $(x_{T+1}^c,y_{T+1}^c)=\textnormal{PB-SSP}(\delta,p,T)$ satisfies $\mathbb{P}\left(\Delta_{\Phi}(x_{T+1}^c,y_{T+1}^c) \leq \epsilon
    \right) \geq 1-p$ and the total computational cost is upper bounded by 
\begin{equation*}
    \mathcal{O}\left(\ln\left(\frac{\ln(\kappa)}{p}\right)\ln(\kappa)\cdot C_{\mathcal{M}}^w\left(\Phi,\frac{\epsilon}{\ln(\kappa)}\right)\right).
    \nonumber
\end{equation*}
\label{corollary: PB-SSP bound}
\end{corollary}
\vspace{-1.0cm}
\proof{proof}
To generate $(x_i^c,y_i^c)$ satisfying (\ref{equ: high prob x and y}) at round $i=0,\dots,T$, we can call $m=\left\lceil18\ln\left(\frac{2T+4}{p}\right)\right\rceil$ times of $\mathcal{M}^w\left(\Phi^{i-1}_x,\delta/27\right)$ and $\mathcal{M}^w(\Phi^{i-1}_y,\delta/27)$ to get two sets of candidate points $X=\{\hat{x}_j\}_{j=1}^m$ and $Y=\{\hat{y}_j\}_{j=1}^m$, respectively. By applying the Markov's inequality and the lower side of ($\ref{equ: unconstrained two-side}$), we can ensure these candidate points satisfy $\mathbb{P}\left[\|\hat{x}_j-x_i^*\|\leq\sqrt{2\delta/9(\mu_x+\lambda_x^{i-1})} \right]\geq2/3$ and $\mathbb{P}\left[\|\hat{y}_j-y_i^*\|\leq\sqrt{2\delta/9(\mu_y+\lambda_y^{i-1})} \right]\geq2/3$. Then we use Algorithm \ref{alg:extract} to extract eligible $(x_i^c,y_i^c)$, and (\ref{equ: high prob x and y}) is guaranteed by applying Lemma \ref{lemma:robust distance estimation}. The way to generate $(x_{T+1}^c,y_{T+1}^c)$ at round $T+1$ is similar, where $\mathcal{M}^w\left(\Phi^{T}_x,\frac{\delta(\mu_x+\lambda_x^T)}{27(L_f+\lambda_x^T)}\right)$ and $\mathcal{M}^w\left(\Phi^{T}_y,\frac{\delta(\mu_y+\lambda_y^T)}{27(L_g+\lambda_y^T)}\right)$ are called, and the upper side of ($\ref{equ: unconstrained two-side}$) is further applied. Noting that both $\frac{L_f+\lambda_x^T}{\mu_x+\lambda_x^T}$ and $\frac{L_g+\lambda_y^T}{\mu_y+\lambda_y^T}$ are less than $2$ according to the above setting, the total sample complexity is less than
\begin{equation*}
    \left\lceil18\ln\left(\frac{2T+4}{p}\right)\right\rceil \left(\sum_{i=0}^T \left(C_{\mathcal{M}}^w(\Phi^{i-1}_x,\delta/27)+C_{\mathcal{M}}^w(\Phi^{i-1}_y,\delta/27)\right)+C_{\mathcal{M}}^w(\Phi^{T}_x,\delta/54)+C_{\mathcal{M}}^w(\Phi^{T}_y,\delta/54)\right).
\end{equation*}
 Moreover, since $\Phi^{i-1}_x$ and $\Phi^{i-1}_y$ are constructed by adding quadratic regularization terms to $\Phi$, it generally holds that
\begin{equation*}
    C_{\mathcal{M}}^w(\Phi^{i-1}_x,\delta) \leq C_{\mathcal{M}}^w(\Phi,\delta) \quad \text{and} \quad
    C_{\mathcal{M}}^w(\Phi^{i-1}_y,\delta) \leq C_{\mathcal{M}}^w(\Phi,\delta), \quad i=0,\dots,T+1. 
\end{equation*}
Hence, the desirable high-probability guarantee and total sample complexity bound are achieved by substituting the above parameter settings.
\Halmos
\endproof
\begin{remark}
    In the analysis of this paper, we choose the base number $\nu$ to be $2$ for ease of deriving the bound. Nevertheless, it can take other values without changing the order of the bound.
\end{remark}

\subsection{Consequences for SAA solutions}
\label{subsec: Consequence of SAA unconstrained}
In this part, we examine the consequences of PB-SSP for a specific SSP oracle based on SAA solutions returned by solving (\ref{equ: saa minmax}). The central question is to determine the total sample size, such that the estimated solution $(\bar{x},\bar{y})$ satisfies the high-probability guarantees as (\ref{equ: high prob guarantee}). Following \cite{zhang2021generalization}, we make one more assumption, and the Lemma \ref{lemma: SAA-oracle unconstrained} characterizing the oracle follows. 
\begin{assumption}
\label{assump: bounded C}
There exists a constant $C$ s.t. $\mathbb{E}_{\xi}[\|\nabla \Phi_{\xi}(x^*,y^*)\|^2]\leq C < +\infty$.
\end{assumption}

\begin{lemma} \textbf{\emph{(\cite[Theorem 4]{zhang2021generalization})}}
\label{lemma: SAA-oracle unconstrained}
Suppose Assumption \ref{assump:SC-SC}, \ref{assump:smoothness} and \ref{assump: bounded C} hold. Given $n$ i.i.d. samples $\{\xi_i\}_{i=1}^n$, the  solution $(\hat{x},\hat{y})$ to the SAA empirical problem \eqref{equ: saa minmax} satisfies 
$$\mathbb{E}[\|\hat{x}-x^*\|^2] \leq \frac{32CL_{xy}^2}{n\mu_x^2\mu_y^2}\qquad \mbox{and} \qquad\mathbb{E}[\|\hat{y}-y^*\|^2] \leq \frac{32CL_{xy}^2}{n\mu_x^2\mu_y^2},$$
where $(x^*,y^*)$ is the saddle point to the true objective function $\Phi$.
\end{lemma}

This result directly bounds the Euclidean distance between the estimated solutions and the optimal solutions. By Markov's inequality, we have $\mathbb{P}\left[\|\hat{x}-x^*\| \leq \sqrt{\frac{96CL_{xy}^2}{n\mu_x^2\mu_y^2}} \right] \geq 2/3$ and $\mathbb{P}\left[\|\hat{y}-y^*\| \leq \sqrt{\frac{96CL_{xy}^2}{n\mu_x^2\mu_y^2}} \right] \geq 2/3$, which enables us to boost the confidence of the SAA oracle for unconstrained SSP problems by incorporating Algorithm \ref{alg:PB-SSP}. Next, we show how we should specify the parameters of the SAA in Algorithm \ref{alg:PB-SSP} so that the confidence of the oracle can be efficiently boosted. For ease of notation, let us define the SAA oracle with proximity terms as follows.

\vspace{0.5cm}
\begin{algorithm}[H]
\SetKwInOut{Return}{Return}
\caption{SAA($n,\lambda_x,\lambda_y,x^c,y^c$)}\label{alg: SAA}
\KwIn{sample count $n \in \mathbb{N}$, amplitude $\lambda_x$, $\lambda_y \geq 0$, center $x^c\in \mathcal{X}$, $y^c\in \mathcal{Y}$.}
Generate i.i.d. samples $\xi_1,\dots,\xi_n \sim \mathcal{P}$ and compute the solution $(\hat{x},\hat{y})$ of
\begin{equation*}
    \min_{x\in \mathcal{X}}\max_{y \in \mathcal{Y}} \frac{1}{n}\sum_{i=1}^n \Phi_{\xi_i}(x,y) + \frac{\lambda_x}{2}\|x-x^c\|^2 -\frac{\lambda_y}{2}\|y-y^c\|^2.
\end{equation*}\\
\Return{$(\hat{x},\hat{y})$} 
\end{algorithm}\vspace{0.5cm}

Now let us show how \eqref{equ: high prob x and y} of the PB-SSP method can be ensured. Let us still adopt previous notations in Algorithm \ref{alg:PB-SSP}. Then we should have $x_{i}^* := \argmin_{x\in\mathcal{X}}\max_{y\in\mathcal{Y}} \,\,  \Phi(x,y) + \frac{\lambda_x^{i-1}}{2}\|x-x^c_{i-1}\|^2$. Now, given an arbitrary batch of $n$ i.i.d. samples $\{\xi_i\}_{i=1}^n$, we set 
$$\hat{x} = \argmin_{x\in\mathcal{X}}\,\,\max_{y\in\mathcal{Y}} \,\, \frac{1}{n}\sum_{j=1}^n \Phi_{\xi_j}(x,y) + \frac{\lambda_x^{i-1}}{2}\|x-x^c_{i-1}\|^2.$$
In other words, we set $(\hat x, \sim) = \text{SAA}(n,\lambda_x^{i-1},0,x_{i-1}^c,\text{null})$. By Lemma \ref{lemma: SAA-oracle unconstrained}, we immediately have 
\begin{eqnarray}
    \label{eqn: saa-sub-1}
    \mathbb{P}\left[\|\hat{x}-x^*_{i}\| \leq \sqrt{\frac{96CL_{xy}^2}{n(\mu_x+\lambda_x^{i-1})^2\mu_y^2}} \right] \geq 2/3.
\end{eqnarray} 
Similarly, with the definition that $y_{i}^* := \argmax_{y\in\mathcal{Y}}\min_{x\in\mathcal{X}} \,\,  \Phi(x,y) - \frac{\lambda_y^{i-1}}{2}\|y-y^c_{i-1}\|^2$ and  we set $(\sim, \hat y) = \text{SAA}(n,0,\lambda_y^{i-1},\text{null},y_{i-1}^c)$, then Lemma \ref{lemma: SAA-oracle unconstrained} also indicates that
\begin{eqnarray}
    \label{eqn: saa-sub-2}
    \mathbb{P}\left[\|\hat{y}-y^*_{i}\| \leq \sqrt{\frac{96CL_{xy}^2}{n\mu_x^2(\mu_y+\lambda_y^i)^2}} \right] \geq 2/3.
\end{eqnarray} 
Consequently, we can apply robust distance estimation to find the high-probability solutions $x_i^c$ and $y_i^c$ that satisfy \eqref{equ: high prob x and y} in Algorithm \ref{alg:PB-SSP}. For ease of notation, we define the following algorithm that combines the SAA oracle and the robust distance estimation.

\begin{algorithm}[htb!]
\SetKwInOut{Return}{Return}
\caption{RobustSAA($n,m,\lambda_x,\lambda_y,x^c,y^c$)}
\label{alg: RobustSAA}
\KwIn{sample count $n \in \mathbb{N}$, trial count $m\in \mathbb{N}$, amplitude $\lambda_x$, $\lambda_y \geq 0$, center $x^c\in \mathcal{X}$, $y^c\in \mathcal{Y}$.}
Let $X=\{\ \}$ and $Y=\{\ \}$ be two empty lists, and set $\rho$ to be Euclidean norm.\\
Call $m$ times SAA($n,\lambda_x,\lambda_y,x^c,y^c$), and add solutions $\{\hat{x}_{j},\hat{y}_{j}\}_{j=1}^m$ into $X$ and $Y$ respectively.\\
Compute $\mathcal{I}_x=\text{Extract}(X,\rho)$, and pick an arbitrary $k_1\in\mathcal{I}_x$.\\
Compute $\mathcal{I}_y=\text{Extract}(Y,\rho)$, and pick an arbitrary $k_2\in\mathcal{I}_y$.\\
\Return{$(\hat{x}_{k_1},\hat{y}_{k_2})$} 
\end{algorithm}

\noindent By Lemma \ref{lemma:robust distance estimation}, \eqref{eqn: saa-sub-1}, and \eqref{eqn: saa-sub-2}, setting 
$(x_i^c,\sim)= \text{RobustSAA}(n_x^{i-1},m,\lambda_{x}^{i-1},0,x_{i-1}^c,\textnormal{null})$ and $(\sim,y_i^c)= \text{RobustSAA}(n_y^{i-1},m,0,\lambda_{y}^{i-1},\textnormal{null}, y_{i-1}^c)$
with $m=\left\lceil18\ln\left(\frac{2T+4}{p}\right)\right\rceil$, $n_x^{i-1}=\left\lceil\frac{432CL_{xy}^2}{(\mu_x+\lambda_x^{i-1})\mu_y^2\delta}\right\rceil$, and $n_y^{i-1}=\left\lceil\frac{432CL_{xy}^2}{\mu_x^2(\mu_y+\lambda_y^{i-1})\delta}\right\rceil$ is sufficient to guarantee \eqref{equ: high prob x and y}.  For \eqref{equ: high prob x and y last round}, the bound can be obtained by finding the high-probability bounds on $\|x^c_{T+1}-x^*_{T+1}\|$ and $\|y^c_{T+1}-y^*_{T+1}\|$ and then applying the upper side of \eqref{equ: unconstrained two-side}. We summarize the discussion as the following algorithm and Theorem \ref{them: Efficiency of BoostSAA} gives the guarantee.\vspace{1cm}

\begin{algorithm}[H]
\SetKwInOut{Return}{Return}
\caption{BoostSAA($\delta, p, T$)}\label{alg:BoostSAA}
\KwIn{$\delta > 0$, $p \in (0,1)$, $T\in \mathbb{N}$.}
Set $\lambda_x^{-1}=\lambda_y^{-1}=0$, $x_{-1}^c=y_{-1}^c=\textnormal{null}$ and $m=\left\lceil18\ln\left(\frac{2T+4}{p}\right)\right\rceil$.\\
\For{$i=0,\dots,T$}{
Set $n_x^{i-1}=\left\lceil\frac{432CL_{xy}^2}{(\mu_x+\lambda_x^{i-1})\mu_y^2\delta}\right\rceil \quad \text{and} \quad
n_y^{i-1}=\left\lceil\frac{432CL_{xy}^2}{\mu_x^2(\mu_y+\lambda_y^{i-1})\delta}\right\rceil$.\\
$x_i^c= \text{RobustSAA}(n_x^{i-1},m,\lambda_{x}^{i-1},0,x_{i-1}^c,\textnormal{null})$,
$y_i^c= \text{RobustSAA}(n_y^{i-1},m,0,\lambda_{y}^{i-1},\textnormal{null}, y_{i-1}^c)$.
}
Set $n_x^T=\left\lceil \frac{L_{xy}^2/\mu_y+L_x+\lambda_x^{T}}{\mu_x+\lambda_x^{T}} \cdot  \frac{432CL_{xy}^2}{(\mu_x+\lambda_x^{T})\mu_y^2\delta}\right\rceil \quad \text{and} \quad
n_y^T=\left\lceil\frac{L_{xy}^2/\mu_x+L_y+\lambda_y^{T}}{\mu_y+\lambda_y^{T}} \cdot \frac{432CL_{xy}^2}{\mu_x^2(\mu_y+\lambda_y^{T})\delta}\right\rceil$.\\
\Return{$
x_{T+1}^c=\text{RobustSAA}\left(n_x^{T},m,\lambda_{x}^{T},0,x_{T}^c,\textnormal{null}\right)$\\
$y_{T+1}^c=\text{RobustSAA}\left(n_y^{T},m,0,\lambda_{y}^{T},\textnormal{null},y_{T}^c\right)
$}
\end{algorithm}\vspace{1cm}

\begin{theorem}[Efficiency of BoostSAA]
Fix a target relative accuracy $\delta>0$, a probability of failure $p \in (0,1)$, and a natural number $T \in \mathbb{N}$. Then with probability at least $1-p$, the point $(x_{T+1}^c,y_{T+1}^c)=\textnormal{BoostSAA}(\delta,p,T)$ satisfies
\begin{equation*}
    \Delta_{\Phi}(x_{T+1}^c,y_{T+1}^c) \leq 
    \delta\left(2+\sum_{i=0}^T\frac{\lambda_x^i}{\mu_x+\lambda_x^{i-1}}+\frac{\lambda_y^i}{\mu_y+\lambda_y^{i-1}}\right).
\end{equation*}
\label{them: Efficiency of BoostSAA}
\end{theorem}
\proof{Proof}
We will verify that Algorithm \ref{alg:BoostSAA} fits the framework of Algorithm \ref{alg:PB-SSP}. Specifically, we check that ($\ref{equ: high prob x and y}$) and (\ref{equ: high prob x and y last round}) are satisfied. During round $i$, $i=0,\dots,T$, Theorem \ref{lemma: SAA-oracle unconstrained} and the definition of $n_x^{i-1}$ guarantee
\begin{equation*}
    \mathbb{P}\left[\|\hat{x}_{j}-x_i^*\| \leq \sqrt{\frac{2\delta}{9(\mu_x+\lambda_x^{i-1})}} \right] \geq 2/3, \quad j=1,\dots,m
\end{equation*}
for all the $\hat{x}_{j}$ in the list $X$ inside the RobustSAA algorithm. Then by a direct application of Lemma \ref{lemma:robust distance estimation}, we deduce that
\begin{equation*}
    \mathbb{P}\left[\|x_{i}^c-x^*_{i}\|\leq \sqrt{\frac{2\delta}{\mu_x+\lambda_x^{i-1}}}\right]\geq 1-\exp\left(-\frac{m}{18}\right) \geq 1-\frac{p}{2T+4},
\end{equation*}
where the second inequality is due to the setting of $m$. In the last round $T+1$, Theorem \ref{lemma: SAA-oracle unconstrained} and the definition of $n_x^{T}$ guarantee
\begin{equation*}
    \mathbb{P}\left[\|\hat{x}_{j}-x_{T+1}^*\| \leq \sqrt{\frac{2\delta}{9(L_{xy}^2/\mu_y+L_x+\lambda_x^{T})}} \right] \geq 2/3, \quad j=1,\dots,m
\end{equation*}
for all the $\hat{x}_{j}$ in the list $X$ inside the RobustSAA algorithm. Lemma 1 guarantees 
\begin{equation*}
    \mathbb{P}\left[\|x_{T+1}^c-x^*_{T+1}\|\leq \sqrt{\frac{2\delta}{L_{xy}^2/\mu_y+L_x+\lambda_x^{T}}}\right]\geq 1-\exp\left(-\frac{m}{18}\right) \geq 1-\frac{p}{2T+4}.
\end{equation*}
Note that, we have $f^T(x_{T+1}^c)-f^T(x^*_{T+1}) \leq \frac{1}{2}\left(L_{xy}^2/\mu_y+L_x+\lambda_x^T\right)\|x_{T+1}^c-x^*_{T+1}\|^2$ due to the Lipschitz smoothness of $f^T(x)$ characterized in (\ref{equ: L of f^i and g^i}). Thus, we deduce that
\begin{equation*}
    \mathbb{P}\left[f^T(x_{T+1}^c)-f^T(x^*_{T+1}) \leq \delta \right]\geq 1-\frac{p}{2T+4}.
\end{equation*}
The other part for $y_i^c$, $i=0,\dots,T+1$ can be proved through a completely parallel way. Thus, both (\ref{equ: high prob x and y}) and (\ref{equ: high prob x and y last round}) are realized in Algorithm \ref{alg:BoostSAA}, and the proof completes.
\Halmos
\endproof

Likewise, the following parameter setting yields the high-probability guarantee.
\begin{corollary}[Efficiency of BoostSAA with geometric decay]
Fix a target accuracy $\epsilon > 0$, and a probability of failure $p \in (0,1)$. Define the algorithm parameters:
\begin{equation*}
\begin{aligned}
    &T=\left \lceil \log_{\nu}\left( \max\left(\frac{L_{xy}^2/\mu_y+L_x}{\mu_x},\frac{L_{xy}^2/\mu_x+L_y}{\mu_y}\right)
    \right) \right \rceil,\\
    &\delta=\frac{\epsilon}{4+4T}, \quad \lambda_x^i=\mu_x \nu^i \quad \text{and} \quad \lambda_y^i=\mu_y \nu^i \quad \text{with} \quad \nu=2 \quad \forall i \in [0,T].
\end{aligned}
\end{equation*}
Then the point $(x_{T+1}^c,y_{T+1}^c)=\textnormal{BoostSAA}(\delta,p,T)$ satisfies
\begin{equation*}
    \mathbb{P}\left(\Delta_{\Phi}(x_{T+1}^c,y_{T+1}^c) \leq \epsilon
    \right) \geq 1-p.
\end{equation*}
\end{corollary}
Moreover, the total number of samples used by the algorithm can be calculated as
\begin{equation*}
\begin{aligned}
    &m \left(\sum_{i=0}^{T} (n_x^{i-1}+n_y^{i-1})
    + n_x^{T} + n_y^{T}\right)\\
    \leq&\frac{432mCL_{xy}^2}{\mu_y^2\delta} \cdot \left(\sum_{i=0}^{T}
    \frac{1}{\mu_x+\lambda_x^{i-1}} + \frac{2}{\mu_x+\lambda_x^{T}}\right)+
    \frac{432mCL_{xy}^2}{\mu_x^2\delta} \cdot \left(\sum_{i=0}^{T}
    \frac{1}{\mu_y+\lambda_y^{i-1}} + \frac{2}{\mu_y+\lambda_y^{T}}\right).
\end{aligned}
\end{equation*}
Noting that
\begin{equation*}
    \sum_{i=0}^{T}
    \frac{1}{\mu_x+\lambda_x^{i-1}} + \frac{2}{\mu_x+\lambda_x^{T}}=
    \frac{1}{\mu_x}+\sum_{i=1}^{T}\frac{1}{\mu_x+\mu_x2^{i-1}}+\frac{2}{\mu_x+\mu_x2^T} \leq \frac{1}{\mu_x}+\sum_{i=1}^{T}\frac{1}{\mu_x2^{i-1}}+\frac{2}{\mu_x2^T}
    \leq \frac{3}{\mu_x},
\end{equation*}
we conclude the total sample complexity is bounded by
\begin{equation*}
    \mathcal{O}\left(\ln(\kappa)\ln\left(\frac{\ln(\kappa)}{p}\right)
    \frac{C\kappa^2}{\mu \epsilon}
    \right).
    \label{equ: Complexity of BoostSAA}
\end{equation*}
The $\kappa^2$ appearing in the last term roots in the SAA oracle characterized by Lemma \ref{lemma: SAA-oracle unconstrained}. Therefore, this dependence cannot be mitigated by our procedure. 

Finally, it is worth noting that if we do not use the inexact proximal point iterations of PB-SSP, we can directly apply the RDE approach to generate $(\bar{x},\bar{y})=\text{RobustSAA}(n,m,0,0,\text{null},\text{null})$, with $n=\left \lceil \frac{432CL_{xy}^2}{\epsilon \mu_x^2\mu_y^2} \left(\frac{L_{xy}^2}{\mu_y}+ \frac{L_{xy}^2}{\mu_x}+L_x+L_y\right) \right \rceil$ and $m=\lceil 18\ln{(2/p)} \rceil$. The high-probability guarantee \eqref{equ: high prob guarantee} is still achieved, but the overall sample complexity will be 
\begin{equation}
    m\cdot n = \mathcal{O}\left(\ln{\left(\frac{1}{p}\right)}\frac{C\kappa^4}{\mu\epsilon}\right),
    \label{equ: complexity of robust estimation}\nonumber
\end{equation}
which is worse than PB-SSP by a $\kappa^2$ factor in terms of sample efficiency.

\section{Extension to constrained problems}
\label{sec: constrained}
We now consider the constrained SSP problem, where both $\mathcal{X}$ and $\mathcal{Y}$ are compact and convex sets, and Assumption \ref{assump: function Lipschitz} is satisfied. 
However, the techniques in Section 3 are not directly extendable in this scenario due to the nonvanishing gradient terms in the two-sided bounds \eqref{equ: constrained_two_side-1} and \eqref{equ: constrained_two_side-2}: $\langle\nabla_x \Phi(x^*,y^*) ,x-x^*\rangle$ and $\langle\nabla_y \Phi(x^*,y^*) ,y-y^*\rangle$. While we can use an SSP oracle to ensure that these quantities are small in expectation, we cannot directly use robust distance estimation to extract candidate points because they are not well-defined pseudometrics. Such a difficulty further prevents us from generating $(x_{T+1}^c,y_{T+1}^c)$ that satisfies (\ref{equ: high prob x and y last round}) in the last round if we still want to follow the framework of PB-SSP. Nevertheless, we can still apply the techniques presented in Section 3 to generate $(x_i^c,y_i^c)$ for round $i=0,\dots,T$, since the additional two quantities are always non-negative due to the optimality conditions, and the lower bound of Equation (\ref{equ: unconstrained two-side}) still holds for constrained SSP problems.

\subsection{Robust estimation for constrained setting}
To overcome the above difficulty, let us consider the following procedure. Assume we have the exact values of $\nabla_x \Phi(x^*,y^*)$ and $\nabla_y \Phi(x^*,y^*)$, we can define two pseudometrics:
\begin{equation*}
    \rho_x(x_1,x_2)=|\langle\nabla_x \Phi(x^*,y^*) ,x_1-x_2\rangle| \quad \text{and} \quad \rho_y(y_1,y_2)=|\langle\nabla_y \Phi(x^*,y^*) ,y_1-y_2\rangle|.
\end{equation*}
Of course, the exact values of $\nabla_x \Phi(x^*,y^*)$ and $\nabla_y \Phi(x^*,y^*)$ are not accessible, so we will replace them with reasonable estimators and define two alternative pseudometrics. We will see later that in order to control the function errors within an acceptable magnitude, it suffices to approximate $\nabla_x \Phi(x^*,y^*)$ and $\nabla_y \Phi(x^*,y^*)$ up to a very loose accuracy, and the extra computation cost is negligible compared with that of calling the SSP oracle.  More formally, we  make one more assumption for the stochastic gradient of $\Phi$.
\begin{assumption}
\label{assump: bounded variance of gradient}
Fix a probability space $(\Omega,\mathcal{F},\mathcal{P})$ and let $G_x: \boldsymbol{R}^{d_x} \times \boldsymbol{R}^{d_y} \times \Omega \rightarrow \boldsymbol{R}^{d_x}$ and $G_y: \boldsymbol{R}^{d_x} \times \boldsymbol{R}^{d_y} \times \Omega \rightarrow \boldsymbol{R}^{d_y}$ be two measurable maps satisfying
\begin{equation*}
\begin{aligned}
    \mathbb{E}_{\xi}G_x(x,y,\xi)=\nabla_x \Phi(x,y), \qquad
    \mathbb{E}_{\xi} \| G_x(x,y,\xi) - \nabla_x \Phi(x,y) \|^2 \leq \sigma_x^2,\\
    \mathbb{E}_{\xi}G_y(x,y,\xi)=\nabla_y \Phi(x,y), \qquad
    \mathbb{E}_{\xi} \| G_y(x,y,\xi) - \nabla_y \Phi(x,y) \|^2 \leq \sigma_y^2,
\end{aligned}
\end{equation*}
where $d_x$ and $d_y$ are the dimensions of $\mathcal{X}$ and $\mathcal{Y}$ respectively.
\end{assumption}
Under this assumption, we can define two gradient oracles $\mathcal{G}_{\sigma_{x}}(\cdot,\cdot,\delta_{G})$ and $\mathcal{G}_{\sigma_{y}}(\cdot,\cdot,\delta_{G})$ as the average of a finite sample of stochastic gradients, i.e., for any $x \in \mathcal{X}$ and $y \in \mathcal{Y}$,
\begin{equation*}
\begin{aligned}
    \mathcal{G}_{\sigma_{x}}(x,y,\delta_{G}):= \frac{1}{n_x}\sum_{i=1}^{n_x} G_x(x,y,\xi_i) \quad \text{where} \quad
    n_x=\left \lceil \frac{3\sigma_x^2}{\delta_{G}^2} \right \rceil,\\
    \mathcal{G}_{\sigma_{y}}(x,y,\delta_{G}):= \frac{1}{n_y}\sum_{i=1}^{n_y} G_y(x,y,\xi_i) \quad \text{where} \quad
    n_y=\left \lceil \frac{3\sigma_y^2}{\delta_{G}^2} \right \rceil.
\end{aligned}
\end{equation*}
Using Markov's inequality, we have
\begin{equation*}
\begin{aligned}
    \mathbb{P}\left[ \left\| \frac{1}{n_x}\sum_{i=1}^{n_x} G_x(x,y,\xi_i) - \nabla_x \Phi(x,y)\right\|^2 \leq \delta_{G}^2 \right] \geq 
    1-\frac{\sigma_x^2/n_x}{\delta_{G}^2} \geq \frac{2}{3},\\
    \mathbb{P}\left[ \left\| \frac{1}{n_y}\sum_{i=1}^{n_y} G_y(x,y,\xi_i) - \nabla_y \Phi(x,y)\right\|^2 \leq \delta_{G}^2 \right] \geq
    1-\frac{\sigma_y^2/n_y}{\delta_{G}^2} \geq \frac{2}{3}.
\end{aligned}
\end{equation*}
Based on the two gradient oracles, we can apply Algorithm \ref{alg:extract} to generate gradient estimators with high-probability guarantees, which is encoded in Algorithm \ref{alg: gradient}.\vspace{0.3cm}

\begin{algorithm}[H]
\SetKwInOut{Return}{Return}
\caption{Gradient$(x,y,\delta_{G},m,\text{Flag})$}\label{alg: gradient}
\KwIn{a point $(x,y)$, $\delta_{G} >0$, trial count $m \in \mathbb{N}$, Flag $\in \{0,1\}$.}
\eIf{\textnormal{Flag} $=1$}{
    Define the map $G:=G_x$, and $\sigma:=\sigma_x$.}
    {Define the map $G:=G_y$, and $\sigma:=\sigma_y$.}
Let $\widetilde{\mathcal{G}}=\{\ \}$ be an empty list and $n=\left \lceil 3\sigma^2/\delta_{G}^2 \right \rceil$.\\
\For{$j=1,\dots,m$}{
Generate i.i.d. samples $\xi_1,\dots,\xi_{n} \sim \mathcal{P}$ and compute
\begin{equation*}
    \mathcal{G}_{\sigma}^j(x,y,\delta_{G})=\frac{1}{n}\sum_{i=1}^{n} G(x,y,\xi_i).
\end{equation*}\\
Add $\mathcal{G}_{\sigma}^j(x,y,\delta_{G})$ into $\widetilde{\mathcal{G}}$.
}
Set $\rho$ to be Euclidean norm and compute $\mathcal{I}_{\mathcal{G}}=\text{Extract}(\widetilde{\mathcal{G}},\rho)$. Pick an arbitrary $k\in\mathcal{I}_{\mathcal{G}}$.\\
\Return{$\mathcal{G}_{\sigma}^k(x,y,\delta_{G})$}
\end{algorithm}\vspace{0.3cm}

Let $\widetilde{\nabla}_x\Phi(x,y)$ and $\widetilde{\nabla}_y\Phi(x,y)$ denote the outputs of Algorithm \ref{alg: gradient} by setting Flag $=1$ and $0$ respectively. A direct application of Lemma \ref{lemma:robust distance estimation} yields
\begin{equation*}
\begin{aligned}
    \mathbb{P}[\|\widetilde{\nabla}_x\Phi(x,y) - \nabla_x \Phi(x,y)\| \leq 3 \delta_{G}] \geq 1-\exp{(-m/18)},\\
    \mathbb{P}[\|\widetilde{\nabla}_y\Phi(x,y) - \nabla_y \Phi(x,y)\| \leq 3 \delta_{G}] \geq 1-\exp{(-m/18)}.
\end{aligned}
\end{equation*}
Also, we can extend the above to robustly estimate the gradient for all perturbed functions. Define $\widetilde{\nabla}_x\Phi_x^i(x,y):=\widetilde{\nabla}_x\Phi(x,y)+\lambda_x^i(x-x_i^c)$ and $\widetilde{\nabla}_y\Phi_y^i(x,y):=\widetilde{\nabla}_y\Phi(x,y)-\lambda_y^i(y-y_i^c)$. The following results hold for $i=-1,0,\dots,T$.
\begin{equation*}
\begin{aligned}
    \mathbb{P}[\|\widetilde{\nabla}_x\Phi_x^i(x,y) - \nabla_x \Phi_x^i(x,y)\| \leq 3 \delta_{G}] \geq 1-\exp{(-m/18)},\\
    \mathbb{P}[\|\widetilde{\nabla}_y\Phi_x^i(x,y) - \nabla_y \Phi_x^i(x,y)\| \leq 3 \delta_{G}] \geq 1-\exp{(-m/18)}.
\end{aligned}
\end{equation*}
We then propose the Algorithm \ref{alg: robust function gap} to robustly estimate the gap $f^{T}(x)-f^T(x^*_{T+1})$ and $g^T(y^*_{T+1})-g^T(y)$, and the theorem follows.

\begin{algorithm}
\SetKwInOut{Return}{Return}
\caption{FunctionGap($\mathcal{M}^w(\cdot,\cdot),\delta,m,\textnormal{Flag}$)}
\label{alg: robust function gap}
\KwIn{Oracle $\mathcal{M}^w(\cdot,\cdot)$, target relative accuracy $\delta > 0$, an odd number $m \in \mathbb{N}$, Flag $\in \{0,1\}$.}
\eIf{\textnormal{Flag} $=1$}{
    Define the objective function $\Phi_T:=\Phi_x^T$, and $\delta_G:=(L_x+\lambda_x^T)\sqrt{\delta/(\mu_x+\lambda_x^T)}$.}
    {Define the objective function $\Phi_T:=\Phi_y^T$, and $\delta_G:=(L_y+\lambda_y^T)\sqrt{\delta/(\mu_y+\lambda_y^T)}$.}
Independently generate $(\hat{x}_1,\hat{y}_1),\dots,(\hat{x}_m,\hat{y}_m)$ by calling $\mathcal{M}^w(\Phi_T,\delta/3)$ such that
\begin{equation*}
    \mathbb{P}[\Delta^w_{\Phi_T}(\hat x_j,\hat y_j) \leq \delta] \geq \frac{2}{3}, \quad \text{for all } j \in [1,m].
\end{equation*}\\
Set $\rho_1=\| \cdot \|$ to be Euclidean norm and compute
\begin{equation*}
\begin{aligned}
    \mathcal{I}_1=\text{Extract}(\{\hat x_j\}_{j=1}^m,\rho_1) \quad \text{and} \quad
    \mathcal{I}_2=\text{Extract}(\{\hat y_j\}_{j=1}^m,\rho_1).
\end{aligned}
\end{equation*}\\
Fix arbitrary $k_1 \in \mathcal{I}_1$, $k_2 \in \mathcal{I}_2$ and set $x_G:=\hat x_{k_1}$, $y_G:=\hat y_{k_2}$.\\
\eIf{\textnormal{Flag} $=1$}{
   Compute $\widetilde{\nabla}\Phi_T(x_G,y_G) = \textnormal{Gradient}\left(x_G,y_G,\delta_G,m,\text{Flag}\right)+\lambda_x^T(x_G-x_T^c)$.\\
   Define the pseudometric $\rho_2(x_1,x_2):=| \langle \widetilde{\nabla}\Phi_T(x_G,y_G), x_1-x_2 \rangle |$ on $\mathcal{X}$.\\
   Compute $\mathcal{I}_3=\text{Extract}(\{\hat x_j\}_{j=1}^m,\rho_2)$, and pick an arbitrary $k_3 \in \mathcal{I}_1 \cap \mathcal{I}_3$.\\
   \Return{$\hat x_{k_3}$}}
    {Compute $\widetilde{\nabla}\Phi_T(x_G,y_G) = \textnormal{Gradient}\left(x_G,y_G,\delta_G,m,\text{Flag}\right)-\lambda_y^T(y_G-y_T^c)$.\\
   Define the pseudometric $\rho_2(y_1,y_2):=| \langle \widetilde{\nabla}\Phi_T(x_G,y_G), y_1-y_2 \rangle |$ on $\mathcal{Y}$.\\
   Compute $\mathcal{I}_3=\text{Extract}(\{\hat y_j\}_{j=1}^m,\rho_2)$, and pick an arbitrary $k_3 \in \mathcal{I}_2 \cap \mathcal{I}_3$.\\
   \Return{$\hat y_{k_3}$}}
\end{algorithm}

\begin{theorem}[Robust function gap estimation]
\label{them: function gap estimation}
With probability at least $1-2\exp(-\frac{m}{18})$, the point $x=\textnormal{FunctionGap}(\mathcal{M}^w(\cdot,\cdot),\delta,m,\textnormal{Flag}=1)$ satisfies the guarantee
\begin{equation*}
\begin{aligned}
    f^T(x)-f^T(x^*_{T+1})\leq \left(3+(18\sqrt{2}+45)\frac{L_x+\lambda_x^T}{\mu_x+\lambda_x^T}+\frac{36L_{xy}}{\sqrt{(\mu_x+\lambda_x^T)\mu_y}}+\frac{9L_{xy}^2}{(\mu_x+\lambda_x^T)\mu_y}
    \right)\delta.
\end{aligned}
\end{equation*}
Likewise, with probability at least $1-2\exp(-\frac{m}{18})$, the point $y=\textnormal{FunctionGap}(\mathcal{M}^w(\cdot,\cdot),\delta,m,\textnormal{Flag}=0)$ satisfies the guarantee

\begin{equation*}
\begin{aligned}
    g^T(y^*_{T+1})-g^T(y)\leq \left(3+(18\sqrt{2}+45)\frac{L_y+\lambda_y^T}{\mu_y+\lambda_y^T}+\frac{36L_{xy}}{\sqrt{\mu_x(\mu_y+\lambda_y^T)}}+\frac{9L_{xy}^2}{\mu_x(\mu_y+\lambda_y^T)}
    \right)\delta.
\end{aligned}
\end{equation*}
\end{theorem}

\proof{Proof}
We first show the detailed proof for the results when Flag $=1$. We have $\Phi_T:=\Phi_x^T$ and denote $\widetilde{\nabla}\Phi_T=\widetilde{\nabla}_x\Phi_x^T$ to avoid confusion.
Define the index set $\mathcal{J}=\{j \in [1,m]:\Delta^w_{\Phi_x^T}(\hat x_j,\hat y_j):=\Phi_x^T(\hat x_j,y^*_{x,T+1})-\Phi_x^T(x^*_{T+1},\hat y_j)\leq \delta \}$ and the event
\begin{equation*}
    E_1:=\left\{|\mathcal{J}|>\frac{m}{2} \right\}.
\end{equation*}
Hoeffding's inequality for Bernoulli random variable guarantees $\mathbb{P}[E_1]\geq 1-\exp(-m/18)$.
Moreover, using the lower bound (\ref{equ: constrained_two_side-1}), we deduce it holds for all $j \in \mathcal{J}$ that
\begin{equation*}
    \|\hat x_j-x^*_{T+1}\|\leq \sqrt{\frac{2\delta}{(\mu_x+\lambda_x^T)}},\quad \|\hat y_j-y^*_{x,T+1}\|\leq \sqrt{\frac{2\delta}{\mu_y}} \quad \text{and} \quad \langle\nabla_x \Phi_x^T(x^*_{T+1},y^*_{x,T+1}) ,\hat x_j-x^*_{T+1}\rangle \leq \delta.
\end{equation*}
Henceforth, suppose that the event $E_1$ occurs. Then Lemma \ref{lemma:robust distance estimation} implies
\begin{equation*}
    \|\hat x_{k_1}-x^*_{T+1}\|\leq 3\sqrt{\frac{2\delta}{(\mu_x+\lambda_x^T)}} \quad \text{for all } k_1 \in \mathcal{I}_1, \quad
    \|\hat y_{k_2}-y^*_{x,T+1}\|\leq 3\sqrt{\frac{2\delta}{\mu_y}} \quad \text{for all } k_2 \in \mathcal{I}_2.
\end{equation*}
Define the event 
\begin{equation*}
    E_2:=\left\{\|\widetilde{\nabla}\Phi_x^T(x_G,y_G) - \nabla_x \Phi_x^T(x_G,y_G)\|\leq 3\delta_G = 3(L_x+\lambda_x^T) \sqrt{\delta/(\mu_x+\lambda_x^T)} \right\},
\end{equation*}
and we know that $\mathbb{P}[E_2]\geq 1-\exp{(-m/18)}$. Suppose that $E_1 \cap E_2$ occurs. Then, we compute
\begin{equation*}
\begin{aligned}
    &\|\widetilde{\nabla}_x\Phi_x^T(x_G,y_G)-\nabla_x\Phi_x^T(x^*_{T+1},y^*_{x,T+1})\|\\
    \leq &\|\widetilde{\nabla}_x\Phi_x^T(x_G,y_G)-\nabla_x\Phi_x^T(x_G,y_G)\|+
    \|\nabla_x\Phi_x^T(x_G,y_G)-\nabla_x\Phi_x^T(x^*_{T+1},y^*_{x,T+1})\| \\
    \leq &3(L_x+\lambda_x^T) \sqrt{\delta/(\mu_x+\lambda_x^T)}+(L_x+\lambda_x^T)\|x_G-x^*_{T+1}\|+L_{xy}\|y_G-y^*_{x,T+1}\| \\
    \leq &3(L_x+\lambda_x^T) \sqrt{\delta/(\mu_x+\lambda_x^T)}+3\sqrt{2}(L_x+\lambda_x^T) \sqrt{\delta/(\mu_x+\lambda_x^T)}+3L_{xy}\sqrt{2\delta/\mu_y}\\
    =& (3+3\sqrt{2})(L_x+\lambda_x^T) \sqrt{\delta/(\mu_x+\lambda_x^T)}+3L_{xy}\sqrt{2\delta/\mu_y}.
\end{aligned}
\end{equation*}
Consequently, for each index $j \in \mathcal{J}$, we successively deduce
\begin{equation*}
\begin{aligned}
    \rho_2(\hat x_j,x^*_{T+1})
    &=|\langle \widetilde{\nabla}_x\Phi_x^T(x_G,y_G), \hat x_j-x^*_{T+1}\rangle|\\
    &\leq \langle\nabla_x\Phi_x^T(x^*_{T+1},y^*_{x,T+1}), \hat x_j-x^*_{T+1}\rangle +
    |\langle \widetilde{\nabla}_x\Phi_x^T(x_G,y_G)-\nabla_x\Phi_x^T(x^*_{T+1},y^*_{x,T+1}), \hat x_j-x^*_{T+1}\rangle|\\
    & \leq \delta + \left((3+3\sqrt{2})(L_x+\lambda_x^T) \sqrt{\delta/(\mu_x+\lambda_x^T)}+3L_{xy}\sqrt{2\delta/\mu_y} \right)\sqrt{2\delta/(\mu_x+\lambda_x^T)}\\
    & = \left(1+(3\sqrt{2}+6)(L_x+\lambda_x^T)/(\mu_x+\lambda_x^T)+6L_{xy}/\sqrt{(\mu_x+\lambda_x^T)\mu_y}
    \right)\delta.
\end{aligned}
\end{equation*}
Therefor, in the event $E_1 \cap E_2$, we conclude
\begin{equation*}
    \rho_2(\hat x_j,x^*_{T+1})\leq \left(3+(9\sqrt{2}+18)(L_x+\lambda_x^T)/(\mu_x+\lambda_x^T)+18L_{xy}/\sqrt{(\mu_x+\lambda_x^T)\mu_y}
    \right)\delta \quad \text{for all } j \in \mathcal{I}_3.
\end{equation*}
Finally, fix an arbitrary index $k_3 \in \mathcal{I}_1 \cap \mathcal{I}_3$. We therefore deduce
\begin{equation*}
\begin{aligned}
    & \langle\nabla_x\Phi_x^T(x^*_{T+1},y^*_{x,T+1}), \hat x_{k_3}-x^*_{T+1}\rangle\\
    \leq &\rho_2(\hat x_{k_3},x^*_{T+1})+|\langle \nabla_x\Phi_x^T(x^*_{T+1},y^*_{x,T+1})- \widetilde{\nabla}_x\Phi_x^T(x_G,y_G), \hat x_{k_3}-x^*_{T+1}\rangle|\\
    \leq &\left(3+(9\sqrt{2}+18)(L_x+\lambda_x^T)/(\mu_x+\lambda_x^T)+18L_{xy}/\sqrt{(\mu_x+\lambda_x^T)\mu_y}
    \right)\delta+\\
    &\left((3+3\sqrt{2})(L_x+\lambda_x^T) \sqrt{\delta/(\mu_x+\lambda_x^T)}+3L_{xy}\sqrt{2\delta/\mu_y} \right)3\sqrt{2\delta/(\mu_x+\lambda_x^T)}\\
    =& \left(3+(18\sqrt{2}+36)(L_x+\lambda_x^T)/(\mu_x+\lambda_x^T)+36L_{xy}/\sqrt{(\mu_x+\lambda_x^T)\mu_y}
    \right)\delta.
\end{aligned}
\end{equation*}
Using the upper bound of (\ref{equ: constrained_two_side-2}), we therefore conclude
\begin{equation*}
\begin{aligned}
    & f^T(\hat x_{k_3})-f^T(x^*_{T+1})\\
    \leq &\frac{1}{2}\left(\frac{L_{xy}^2}{\mu_y}+L_x+\lambda_x^T\right)\|\hat x_{k_3}-x^*_{T+1}\|^2+\langle\nabla_x\Phi_x^T(x^*_{T+1},y^*_{x,T+1}), \hat x_{k_3}-x^*_{T+1}\rangle \\
    \leq &\left(3+ (18\sqrt{2}+45)(L_x+\lambda_x^T)/(\mu_x+\lambda_x^T)+36L_{xy}/\sqrt{(\mu_x+\lambda_x^T)\mu_y}+9L_{xy}^2/((\mu_x+\lambda_x^T)\mu_y)
    \right)\delta.
\end{aligned}
\end{equation*}
Noting $\mathbb{P}[E_1 \cap E_2] \geq 1-2\exp(-m/18)$, the proof for Flag $=1$ completes. The results of the second part can be proved in a completely parallel way.
\Halmos
\endproof

The intuition behind Algorithm \ref{alg: robust function gap} is that, the returned $\hat x_{k_3}$ with $k_3\in \mathcal{I}_1 \cap \mathcal{I}_3$ simultaneously achieves low values of $\|\hat x_{k_3}-x^*_{T+1}\|$ and $\langle\nabla_x\Phi(x^*_{T+1},y^*_{x,T+1}), \hat x_{k_3}-x^*_{T+1}\rangle$ with high probability. Recall that Algorithm \ref{alg:extract} produces at least $m/2$ points, so $m$ must be an odd number to ensure $\mathcal{I}_1 \cap \mathcal{I}_3$ is not empty. Similar arguments hold for the returned $y_k$. With Algorithm \ref{alg: robust function gap} on hand, we can generate $(x_{T+1}^c,y_{T+1}^c)$ satisfying (\ref{equ: high prob x and y last round}).
Besides, note that Algorithm \ref{alg: robust function gap} can also be applied to robustly estimate $f(x)-f(x^*)$ and $g(y^*)-g(y)$. To make it happen, we can simply change the index $T$ to $-1$, set $\lambda_x^{-1}=\lambda_y^{-1}=0$ such that $x_0^*=x^*$, $y_0^*=y^*$, and Theorem \ref{them: function gap estimation} adapts accordingly.

\subsection{Consequences for SAA solutions}\label{subsec: constrained SAA}
In this part, we examine the results for a specific SSP oracle based on SAA solutions. We will still follow the general framework of PB-SSP outlined in Section \ref{subsec: PB-SSP}, whereas Algorithm \ref{alg: robust function gap} will be in place to generate $(x_{T+1}^c,y_{T+1}^c)$ in the last round. To formally characterize this SAA oracle for constrained SSP problems, we provide the following lemma.

\begin{lemma}[{\cite[Theorem 1 and Lemma 2]{zhang2021generalization}}]
\label{lemma: SAA bounded oracle}
Fix an i.i.d. sample $\xi_1,\xi_2,\dots,\xi_n$. Under Assumption \ref{assump:SC-SC} and \ref{assump: function Lipschitz}, the solution $(\hat{x},\hat{y})$ to the SAA problem $\hat{\Phi}_n$ defined by \eqref{equ: saa minmax} satisfies the bound:
\begin{equation*}
    \mathbb{E}[\Delta^w_{\Phi}(\hat{x},\hat{y})] \leq \frac{2}{n}\left(\frac{\ell_x^2}{\mu_x}+\frac{\ell_y^2}{\mu_y}
    \right).
\end{equation*}
Moreover, denote $(\hat{x},\hat{y}_{x})$ and $(\hat{x}_{y},\hat{y})$ the solutions to the SAA problems adding regularization terms $\hat{\Phi}_{n}+\frac{\lambda^i_x}{2}\|x-x_i^c\|^2$ and $\hat{\Phi}_{n}-\frac{\lambda^i_y}{2}\|y-y_i^c\|^2$ respectively. They are bounded by
\begin{equation*}
\begin{aligned}
    &\mathbb{E}[\Delta^w_{\Phi_x^i}(\hat{x},\hat{y}_{x})] \leq \frac{2}{n}\left(\frac{\ell_x^2}{\mu_x+\lambda_x^i}+\frac{\ell_y^2}{\mu_y}\right),\\
    &\mathbb{E}[\Delta^w_{\Phi_y^i}(\hat{x}_{y},\hat{y})] \leq \frac{2}{n}\left(\frac{\ell_x^2}{\mu_x}+\frac{\ell_y^2}{\mu_y+\lambda_y^i}\right).\\
\end{aligned}
\end{equation*}
\end{lemma}
Note that $\ell_x$ and $\ell_y$ originate from $\Phi$, and the Lipschitz constant of the regularization terms will not contribute to the above bounds.

The above oracle can bound the weaker variant of the duality gap. Use Algorithm \ref{alg: SAA} and \ref{alg: RobustSAA} in Section 3 without any change, and consider the returned solution $(\hat x,\sim)=\text{SAA}(n,\lambda_x^{i-1},0,x_{i-1}^c,\text{null})$. By Markov's inequality and applying the lower bound \eqref{equ: constrained_two_side-1}, Lemma \ref{lemma: SAA bounded oracle} implies
\begin{eqnarray}
    \label{eqn: saa-cons-x}
    \mathbb{P}\left[\|\hat{x}-x^*_{i}\| \leq \sqrt{\frac{12}{n}\left(\frac{\ell_x^2}{(\mu_x+\lambda_x^{i-1})^2}+\frac{\ell_y^2}{(\mu_x+\lambda_x^{i-1})\mu_y}\right)} \right] \geq 2/3.
\end{eqnarray}
Similarly, setting $(\sim,\hat y)=\text{SAA}(n,0,\lambda_y^{i-1},\text{null},y_{i-1}^c)$, we deduce
\begin{eqnarray}
    \label{eqn: saa-cons-y}
    \mathbb{P}\left[\|\hat{y}-y^*_{i}\| \leq \sqrt{\frac{12}{n}\left(\frac{\ell_x^2}{\mu_x(\mu_y+\lambda_y^{i-1})}+\frac{\ell_y^2}{(\mu_y+\lambda_y^{i-1})^2}\right)} \right] \geq 2/3.
\end{eqnarray}
Consequently, we can leverage the robust distance estimation to generate $(x_i^c,y_i^c)$ satisfying \eqref{equ: high prob x and y} in Algorithm \ref{alg:PB-SSP}. That is, setting 
$(x_i^c,\sim)= \text{RobustSAA}(n_x^{i-1},m,\lambda_{x}^{i-1},0,x_{i-1}^c,\textnormal{null})$ and $(\sim,y_i^c)= \text{RobustSAA}(n_y^{i-1},m,0,\lambda_{y}^{i-1},\textnormal{null}, y_{i-1}^c)$ with $m=\left\lceil18\ln\left(\frac{2T+4}{p}\right)\right\rceil$, $n_x^{i-1}=\left\lceil\frac{54}{\delta}\left(\frac{\ell_x^2}{\mu_x+\lambda_x^{i-1}}+\frac{\ell_y^2}{\mu_y}\right)\right\rceil$, and $n_y^{i-1}=\left\lceil\frac{54}{\delta}\left(\frac{\ell_x^2}{\mu_x}+\frac{\ell_y^2}{\mu_y+\lambda_y^{i-1}}\right)\right\rceil$. Then, Lemma \ref{lemma:robust distance estimation}, and the above two inequalities \eqref{eqn: saa-cons-x} and \eqref{eqn: saa-cons-y} guarantee \eqref{equ: high prob x and y}.

To generate $(x_{T+1}^c,y_{T+1}^c)$ satisfying \eqref{equ: high prob x and y last round}, we need to use Algorithm \ref{alg: robust function gap}.
Recall the definition of SSP oracles $\mathcal{M}$ in Section \ref{subsection: main contribution}, we define the following two operations of the above SAA oracle when applied to $\Phi_x^T$ and $\Phi_y^T$.
\begin{equation}
\begin{aligned}
    &\mathcal{M}^w(\Phi_x^T,\delta):=\text{SAA}(n,\lambda_x^T,0,x_T^c,\textnormal{null}) \quad \text{with} \quad C_{\mathcal{M}}^w(\Phi_x^T,\delta):=n=\left\lceil\frac{2}{\delta}\left(\frac{\ell_x^2}{\mu_x+\lambda_x^T}+\frac{\ell_y^2}{\mu_y}\right)\right\rceil,\\
    &\mathcal{M}^w(\Phi_y^T,\delta):=\text{SAA}(n,0,\lambda_y^T,\textnormal{null},y_T^c) \quad \text{with} \quad C_{\mathcal{M}}^w(\Phi_y^T,\delta):=n=\left\lceil\frac{2}{\delta}\left(\frac{\ell_x^2}{\mu_x}+\frac{\ell_y^2}{\mu_y+\lambda_y^T}\right)\right\rceil.
    \label{equ:saa oracle operation and cost}
\end{aligned}
\end{equation}
To estimate the gradient, we simply declare
\begin{equation*}
    G_x(x,y,\xi):=\nabla_x\Phi_{\xi}(x,y) \quad \text{and} \quad
    G_y(x,y,\xi):=\nabla_y\Phi_{\xi}(x,y).
\end{equation*}
Then we can upper-bound the variance by the second moment
\begin{equation*}
\begin{aligned}
    &\mathbb{E}_{\xi}\|G_x(x,y,\xi)-\nabla_x\Phi(x,y)\|^2\leq
    2\left(\mathbb{E}_{\xi}\|\nabla_x\Phi_{\xi}(x,y)\|^2+\mathbb{E}_{\xi}\|\nabla_x\Phi(x,y)\|^2\right)\leq 4\ell_x^2,\\
    &\mathbb{E}_{\xi}\|G_y(x,y,\xi)-\nabla_y\Phi(x,y)\|^2\leq
    2\left(\mathbb{E}_{\xi}\|\nabla_y\Phi_{\xi}(x,y)\|^2+\mathbb{E}_{\xi}\|\nabla_y\Phi(x,y)\|^2\right)\leq 4\ell_y^2.
\end{aligned}
\end{equation*}
In addition, we define the following two quantities for ease of notation.
\begin{equation*}
\begin{aligned}
    &M_x^T:=3+(18\sqrt{2}+45)(L_x+\lambda_x^T)/(\mu_x+\lambda_x^T)+36L_{xy}/\sqrt{(\mu_x+\lambda_x^T)\mu_y}+9L_{xy}^2/\left((\mu_x+\lambda_x^T)\mu_y\right),\\
    &M_y^T:=3+(18\sqrt{2}+45)(L_y+\lambda_y^T)/(\mu_y+\lambda_y^T)+36L_{xy}/\sqrt{\mu_x(\mu_y+\lambda_y^T)}+9L_{xy}^2/\left(\mu_x(\mu_y+\lambda_y^T)\right).
\end{aligned}
\end{equation*}
By letting $\lambda_x^T$ and $\lambda_y^T$ sufficient large, $M_x^T$ and $M_y^T$ will reduce to $\mathcal{O}(1)$.
Next, we can call Algorithm \ref{alg: robust function gap} to independently generate 
$x_{T+1}^c=\textnormal{FunctionGap}(\mathcal{M}^w(\cdot,\cdot),\delta/M_x^T,m,\textnormal{Flag}=1)$ and $y_{T+1}^c=\textnormal{FunctionGap}(\mathcal{M}^w(\cdot,\cdot),\delta/M_y^T,m,\textnormal{Flag}=0)$ with 
$\mathcal{M}^w(\cdot,\cdot)$ defined in \eqref{equ:saa oracle operation and cost} and $m=\left\lceil18\ln\left(\frac{2T+4}{p}\right)\right\rceil$. Theorem \ref{them: function gap estimation} guarantees \eqref{equ: high prob x and y last round}, with a slight deviation that the right-hand sides of the two inequalities become $1-\frac{2p}{2T+4}$. The reason is that Algorithm \ref{alg: robust function gap} and Theorem \ref{them: function gap estimation} rely on the occurrence of two high-probability events, while robust distance estimation only depends on one. We can easily fix this by modifying $m$ to be $\left\lceil18\ln\left(\frac{2T+6}{p}\right)\right\rceil$ for all rounds. This way, the probability of all events occurring is still lower bounded by $1-p$, and the order of computation cost remains unaffected. Besides, recall that $m$ must be set as an odd number. We summarize the above discussion into the Algorithm \ref{alg:BoostSAAC}. The following theorem and its corollary are immediate consequences of Theorem \ref{them: PB-SSP} and \ref{them: function gap estimation}.\vspace{0.3cm}
\begin{algorithm}[H]
\SetKwInOut{Return}{Return}
\caption{BoostSAA-C($\delta,p,T$)}\label{alg:BoostSAAC}
\KwIn{$\delta > 0$, $p\in(0,1)$, $T \in \mathbb{N}$}
Set $\lambda_x^{-1}=\lambda_y^{-1}=0$, $x_{-1}^c=y_{-1}^c=0$ and $m=\left\lceil18\ln\left(\frac{2T+6}{p}\right)\right\rceil$ or $\left\lceil18\ln\left(\frac{2T+6}{p}\right)\right\rceil+1$.\\
\For{$i=0,\dots,T$}{
Set $n_x^{i-1}=\left\lceil\frac{54}{\delta}\left(\frac{\ell_x^2}{\mu_x+\lambda_x^{i-1}}+\frac{\ell_y^2}{\mu_y}\right)\right\rceil \quad \text{and} \quad
n_y^{i-1}=\left\lceil\frac{54}{\delta}\left(\frac{\ell_x^2}{\mu_x}+\frac{\ell_y^2}{\mu_y+\lambda_y^{i-1}}\right)\right\rceil$.\\

$x_i^c= \text{RobustSAA}(n_x^{i-1},m,\lambda_{x}^{i-1},0,x_{i-1}^c,\textnormal{null})$,
$y_i^c= \text{RobustSAA}(n_y^{i-1},m,0,\lambda_{y}^{i-1},\textnormal{null},y_{i-1}^c)$.}

\Return{$x_{T+1}^c=\text{FunctionGap}\left(\mathcal{M}^w(\cdot,\cdot),\delta/M_x^T,m, \text{Flag} =1\right)$\\ $y_{T+1}^c=\text{FunctionGap}\left(\mathcal{M}^w(\cdot,\cdot),\delta/M_y^T,m,\text{Flag} =0\right)$}
\end{algorithm}\vspace{0.3cm}


\begin{theorem}[Efficiency of BoostSAA-C]
Fix a target relative accuracy $\delta>0$, a probability of failure $p\in(0,1)$, and natural numbers $T \in \mathbb{N}$. Then with probability at least $1-p$, the point $(x_{T+1}^c,y_{T+1}^c)=\textnormal{BoostSAA-C}(\delta,p,T)$ satisfies
\begin{equation*}
    \Delta_{\Phi}(x_{T+1}^c,y_{T+1}^c) \leq 
    \delta\left(2+\sum_{i=0}^T\frac{\lambda_x^i}{\mu_x+\lambda_x^{i-1}}+\frac{\lambda_y^i}{\mu_y+\lambda_y^{i-1}}\right).
\end{equation*}
\end{theorem}

\begin{corollary}[Efficiency of BoostSAA-C with geometric decay]
Fix a target accuracy $\epsilon > 0$, and a probability of failure $p \in (0,1)$. Define the algorithm parameters:
\begin{equation*}
\begin{aligned}
    &T=\left \lceil \log_{\nu}\left( \max\left(\frac{L_x}{\mu_x},\frac{L_y}{\mu_y},\frac{L_{xy}^2}{\mu_x\mu_y}\right)
    \right) \right \rceil,\\
    &\delta=\frac{\epsilon}{4+4T}, \quad \lambda_x^i=\mu_x \nu^i \quad \text{and} \quad \lambda_y^i=\mu_y \nu^i \quad \text{with} \quad \nu=2 \quad \forall i \in [0,T].
\end{aligned}
\end{equation*}
Then the point $(x_{T+1}^c,y_{T+1}^c)=\textnormal{BoostSAA-C}(\delta,p,T)$ satisfies
\begin{equation*}
    \mathbb{P}\left(\Delta_{\Phi}(x_{T+1}^c,y_{T+1}^c) \leq \epsilon
    \right) \geq 1-p.
\end{equation*}
\end{corollary}

Given the above parameter setting, we have $(L_x+\lambda_x^T)/(\mu_x+\lambda_x^T)\leq2$, $(L_y+\lambda_y^T)/(\mu_y+\lambda_y^T)\leq2$,$L_{xy}^2/\left((\mu_x+\lambda_x^T)\mu_y\right)\leq 1$, and $L_{xy}^2/\left(\mu_x(\mu_y+\lambda_y^T)\right)\leq 1$. Hence, 
$M_x^T$ and $M_y^T$ are both upper bounded by $138+36\sqrt{2}$, and the total number of samples used by Algorithm \ref{alg:BoostSAAC} can be calculated as
\begin{equation*}
\begin{aligned}
    &m\left(\sum_{i=0}^T(n_x^{i-1}+n_y^{i-1})+C_{\mathcal{M}}^w\left(\Phi_x^T,\frac{\delta}{3M_x^T}\right)+C_{\mathcal{M}}^w\left(\Phi_y^T,\frac{\delta}{3M_y^T}\right)\right)+\\
    &m\left(\frac{3M_x^T(\mu_x+\lambda_x^T)\sigma_x^2}{(L_x+\lambda_x^T)^2\delta}
    +\frac{3M_y^T(\mu_y+\lambda_y^T)\sigma_y^2}{(L_y+\lambda_y^T)^2\delta}
    \right)\\
    = & m\sum_{i=0}^T\left(\frac{54}{\delta}\left(\frac{\ell_x^2}{\mu_x+\lambda_x^{i-1}}+\frac{\ell_y^2}{\mu_y}+\frac{\ell_x^2}{\mu_x}+\frac{\ell_y^2}{\mu_y+\lambda_y^{i-1}}\right)
    \right)+\\
    &m\left(\frac{6M_x^T}{\delta}\left(\frac{\ell_x^2}{\mu_x+\lambda_x^T}+\frac{\ell_y^2}{\mu_y}\right)
    +\frac{6M_y^T}{\delta}\left(\frac{\ell_x^2}{\mu_x}+\frac{\ell_y^2}{\mu_y+\lambda_y^T}\right)+\frac{3M_x^T(\mu_x+\lambda_x^T)\sigma_x^2}{(L_x+\lambda_x^T)^2\delta}
    +\frac{3M_y^T(\mu_y+\lambda_y^T)\sigma_y^2}{(L_y+\lambda_y^T)^2\delta}\right)\\
    \leq & \frac{m}{\delta} \left(\sum_{i=0}^T\left(\frac{108\ell_x^2}{\mu_x}+\frac{108\ell_y^2}{\mu_y}\right)+\frac{12M_x^T\ell_x^2}{\mu_x}+\frac{12M_y^T\ell_y^2}{\mu_y}+\frac{6M_x^T\ell_x^2}{L_x}+\frac{6M_y^T\ell_y^2}{L_x}\right)\\
    \leq & \frac{m}{\delta} \left(\frac{(108T+12M_x^T)\ell_x^2}{\mu_x}+\frac{(108T+12M_y^T)\ell_y^2}{\mu_y}+\frac{6M_x^T\ell_x^2}{L_x}+\frac{6M_y^T\ell_y^2}{L_x}\right).
\end{aligned}
\end{equation*}
In the first row, the two terms of $C_{\mathcal{M}}^w(\cdot,\cdot)$ as defined in \eqref{equ:saa oracle operation and cost} correspond to the cost of calling the SAA oracle in the last round, and the last two terms measure the cost to execute the two gradient oracles.
After substituting the parameters, we conclude the total sample complexity is bounded by
\begin{equation}
    \mathcal{O}\left(\ln^2(\kappa)\ln\left(\frac{\ln(\kappa)}{p}\right)
    \frac{\ell^2}{\mu \epsilon}
    \right).
    \label{equ: complexity of BoostSAAC}
\end{equation}

As before, we can also directly use Algorithm \ref{alg: robust function gap} to generate $(\bar{x},\bar{y})$ with high-probability guarantees, which is the RDE approach for constrained settings. Denote
\begin{equation*}
\begin{aligned}
    &M_x:=3+(18\sqrt{2}+45)L_x/\mu_x+36L_{xy}/\sqrt{\mu_x\mu_y}+9L_{xy}^2/\left(\mu_x\mu_y\right),\\
    &M_y:=3+(18\sqrt{2}+45)L_y/\mu_y+36L_{xy}/\sqrt{\mu_x\mu_y}+9L_{xy}^2/\left(\mu_x\mu_y\right),
\end{aligned}
\end{equation*}
both of which have the order of $\mathcal{O}(\kappa^2)$. Let $m=\lceil 18\ln{(4/p)} \rceil$ or $m=\lceil 18\ln{(4/p)} \rceil$+1, $\delta=\epsilon/(M_x+M_y)$. We can then independently generate
$\bar{x}=\text{FunctionGap}(\mathcal{M}^w(\cdot,\cdot),\delta,m,\text{Flag} =1)$ and $\bar{y}=\text{FunctionGap}(\mathcal{M}^w(\cdot,\cdot),\delta,m,\text{Flag} =0)$ by letting all index $T$ be $-1$ in Algorithm \ref{alg: robust function gap},
while the high-probability guarantee of (\ref{equ: high prob guarantee}) is ensured by Theorem \ref{them: function gap estimation}.
Noting the two executions of Algorithm \ref{alg: robust function gap} can share the effort for calling the oracle, we calculate the number of samples used to be
\begin{equation*}
    m\cdot C_{\mathcal{M}}^w\left(\Phi,\frac{\epsilon}{3(M_x+M_y)}\right) + m \cdot (M_x+M_y)\cdot \left(\frac{3\sigma_x^2\mu_x}{L^2_x\epsilon}+\frac{3\sigma_y^2\mu_y}{L^2_y\epsilon}\right),
\end{equation*}
which is on the order of 
\begin{equation}
    \mathcal{O}\left(\ln{\left(\frac{1}{p}\right)}\frac{\ell^2\kappa^2}{\mu\epsilon}\right).
    \label{equ: complexity of robust estimation of constrained}\nonumber
\end{equation}
Similar to the results in Section \ref{sec: unconstrained}, this RDE approach is worse than our PB-SSP due to the extra $\kappa^2$ factor of the sample complexity.

\subsection{Extension to general convex and concave problems}
\label{sec: constrained C-C}
Now we extend our analysis of the above SSA oracle to C-C SSP problems, i.e., $\mu_x=\mu_y=0$. To make the above results applicable, we consider an alternative objective function by adding regularization terms to $\Phi$, given by
\begin{equation*}
    \Phi_{\alpha}(x,y):= \Phi(x,y) + \frac{\alpha_x}{2}\|x-x'\|^2- \frac{\alpha_y}{2}\|y-y'\|^2,
\end{equation*}
where $(x',y')$ is an arbitrary point in $\mathcal{X}\times\mathcal{Y}$. Obviously, $\Phi_{\alpha}$ is SC-SC modulus $(\alpha_x,\alpha_y)$.
Denote the diameter of $\mathcal{X}$ and $\mathcal{Y}$ by $D_x$ and $D_y$ respectively, and define $D^2:=D_x^2+D_y^2$. We have the following inequality about duality gaps for any returned solution $(\bar{x},\bar{y})$, given by
\begin{equation*}
\begin{aligned}
    &\Delta_{\Phi_{\alpha}}(\bar{x},\bar{y})\\
    =&\left(\max_{y\in \mathcal{Y}}\Phi(\bar{x},y)+ \frac{\alpha_x}{2}\|\bar{x}-x'\|^2- \frac{\alpha_y}{2}\|y-y'\|^2\right)-
    \left(\min_{x\in \mathcal{X}}\Phi(x,\bar{y})+ \frac{\alpha_x}{2} \|x-x'\|^2- \frac{\alpha_y}{2}\|\bar{y}-y'\|^2\right)\\
    \geq & \left(\max_{y\in \mathcal{Y}}\Phi(\bar{x},y)- \frac{\alpha_y}{2}\|y-y'\|^2\right)-
    \left(\min_{x\in \mathcal{X}}\Phi(x,\bar{y})+ \frac{\alpha_x}{2} \|x-x'\|^2\right)\\
    \geq & \max_{y\in \mathcal{Y}}\Phi(\bar{x},y) - \min_{x\in \mathcal{X}}\Phi(x,\bar{y}) - \frac{(\alpha_xD_x^2+\alpha_yD_y^2)}{2}\\
    = & \Delta_{\Phi}(\bar{x},\bar{y}) - \frac{(\alpha_xD_x^2+\alpha_yD_y^2)}{2}.
\end{aligned}
\end{equation*}
Hence, we can set an alternative accuracy $\epsilon'>0$ and apply Algorithm \ref{alg:BoostSAAC} to $\Phi_{\alpha}$ instead. The returned solution $(\bar{x},\bar{y})$ satisfies
\begin{equation*}
    \mathbb{P}\left(\Delta_{\Phi}(\bar{x},\bar{y}) \leq \epsilon' +(\alpha_x D_x^2+\alpha_y D_y^2)/2
    \right) \geq 1-p.
\end{equation*}
After substituting some parameters in (\ref{equ: complexity of BoostSAAC}), we deduce the sample complexity is bounded by
\begin{equation*}
    \mathcal{O}\left(\ln^2(L/\alpha)\ln\left(\frac{\ln(L/\alpha)}{p}\right)
    \frac{\ell^2}{\alpha \epsilon'}
    \right),
\end{equation*}
where $\alpha:=\min(\alpha_x,\alpha_y)$.
Recalling the remark below Lemma \ref{lemma: SAA bounded oracle}, $\ell$ appearing in the above bound is not affected by the regularization terms. Finally, the following corollary illustrates how to properly set the regularization terms and the alternative accuracy, and establishes the high-probability guarantees for C-C SSP problems.
\begin{corollary}[Efficiency of BoostSAA-C to C-C problems]
Fix a target accuracy $\epsilon > 0$, and a probability of failure $p \in (0,1)$. Set an alternative accuracy $\epsilon'=\epsilon/2$, the regularization parameter $\alpha_x=\epsilon/2D_x^2$ and $\alpha_y=\epsilon/2D_y^2$.
The point $(\bar{x},\bar{y})$ returned by applying Algorithm 
\ref{alg:BoostSAAC} to $\Phi_{\alpha}$ satisfies
\begin{equation*}
    \mathbb{P}\left(\Delta_{\Phi}(\bar{x},\bar{y}) \leq \epsilon
    \right) \geq 1-p,
\end{equation*}
and the total number of samples is bounded by
\begin{equation*}
    \mathcal{O}\left(\ln^2(LD^2/\epsilon)\ln\left(\frac{\ln(LD^2/\epsilon)}{p}\right)
    \frac{\ell^2D^2}{\epsilon^2}
    \right).
\end{equation*}
\end{corollary}
Meanwhile, directly applying Algorithm \ref{alg: robust function gap} to $\Phi_{\alpha}$ and properly choosing parameters can also provide high-probability guarantees, while the sample complexity is on the order of 
\begin{equation*}
    \mathcal{O}\left(\ln\left(\frac{1}{p}\right)
    \frac{\ell^2L^2D^6}{\epsilon^4}
    \right).
\end{equation*}

\section{Compatibility with first-order oracles}
The previous discussion is mainly based on the SAA oracle, which works by first constructing an empirical SSP problem and then solving it by arbitrary algorithms such as extra gradient method. However, as PB-SSP is a general framework that may accommodate arbitrary qualified oracles, in this section, we would like to illustrate its compatibility with stochastic first-order oracles for SC-SC unconstrained and constrained SSP problems, respectively.
\subsection{Unconstrained SSP problems}
For unconstrained SSP problems, we consider the multistage stochastic optimistic gradient descent ascent (MOGDA) method by Fallah et al. \cite{fallah2020optimal} as our oracle because it shares identical assumptions with us. The MOGDA method can provide in expectation guarantees for the squared distance to the saddle point, and we summarize it in the following lemma.
\begin{lemma} \textbf{\emph{(\cite[Corollary V.3.]{fallah2020optimal})}}
\label{lemma: SA-oracle unconstrained}
Suppose Assumption \ref{assump:SC-SC}, \ref{assump:smoothness} and \ref{assump: bounded variance of gradient} hold. Denote $\sigma^2:=\max\{\sigma_x^2,\sigma_y^2\}$ and $\Delta_{\mathrm{in}}:=\|x_0-x^*\|^2+\|y_0-y^*\|^2$. Let $\left(x_n,y_n\right)$ be the solution returned by the MOGDA method after taking $n$ samples, then it  satisfies \vspace{0.5cm}
$$\mathbb{E}\Big[\|x_n-x^*\|^2+\|y_n-y^*\|^2\Big] \leq \mathcal{O}\left(\exp\left(-\frac{\Theta(n)}{\kappa}\right)\Delta_{\mathrm{in}}+\frac{\sigma^2}{n\mu^2}\right)$$
where $(x_0,y_0)$ is the initial point and $(x^*,y^*)$ is the saddle point to the true objective function $\Phi$.
\end{lemma}
The above lemma implies a sample complexity of $\mathcal{O}\left(\kappa\ln\left(\frac{\Delta_{\mathrm{in}}}{\epsilon}\right)+\frac{\sigma^2}{\mu^2\epsilon}\right)$ to guarantee an $\epsilon$-solution in the sense of squared distance. Similar to our discussion in Section \ref{subsec: Consequence of SAA unconstrained}, we can leverage the framework of Algorithm \ref{alg:PB-SSP} to equip MOGDA with the high-probability guarantee for the duality gap. A minor issue is that we need to provide an initial point to the oracle when solving the sequence of subproblems. Naturally, we choose to adopt the solution obtained from the previous iteration as the initial point for the next subproblem. Denote $\text{MOGDA}(\delta,\Phi,x_0,y_0)$ as the operation to use the MOGDA algorithm to generate a solution with $\delta$ small squared distance in expectation starting from the initial point $(x_0,y_0)$. We summarize the procedure to equip MOGDA with high-probability guarantees in the following algorithms, and the corollary is a direct application of Corollary \ref{corollary: PB-SSP bound}.

\vspace{2.5cm}
\begin{algorithm}[H]
\SetKwInOut{Return}{Return}
\caption{RobustMOGDA($\delta,\Phi,x_0,y_0,m$)}
\label{alg: RobustMOGDA}
\KwIn{$\delta > 0$, $\Phi$, $x_0\in \mathcal{X}$, $y_0\in \mathcal{Y}$, $m\in \mathbb{N}$.}
Let $Z=\{\ \}$ be an empty lists, and set $\rho$ to be Euclidean norm.\\
Call $m$ times MOGDA($\delta,\Phi,x_0,y_0$), and add solutions $\{(\hat{x}_{j},\hat{y}_{j})\}_{j=1}^m$ into $Z$.\\
Compute $\mathcal{I}=\text{Extract}(Z,\rho)$, and pick an arbitrary $k\in\mathcal{I}$.\\
\Return{$(\hat{x}_{k},\hat{y}_{k})$} 
\end{algorithm}\vspace{0.3cm}

\begin{algorithm}[H]
\SetKwInOut{Return}{Return}
\caption{BoostMOGDA($\delta, p, T, x_0, y_0$)}\label{alg:BoostMOGDA}
\KwIn{$\delta > 0$, $p \in (0,1)$, $T\in \mathbb{N}$, $x_0\in \mathcal{X}$, $y_0\in \mathcal{Y}$.}
Set $\lambda_x^{-1}=\lambda_y^{-1}=0$, $x_{-1}^c=y_{-1}^c=\textnormal{null}$, $x_{0,x}^{-1}=x_{0,y}^{-1}=x_0$, $y_{0,x}^{-1}=y_{0,y}^{-1}=y_0$ and $m=\left\lceil18\ln\left(\frac{2T+4}{p}\right)\right\rceil$.\\
\For{$i=0,\dots,T$}{
Set $\delta_x^{i-1}=\frac{2\delta}{27(\mu_x+\lambda_x^{i-1})} \quad \text{and} \quad
\delta_y^{i-1}=\frac{2\delta}{27(\mu_y+\lambda_y^{i-1})}$.\\
$(x_i^c,y_{x,i}^c)= \text{RobustMOGDA}\left(\delta_x^{i-1},\Phi_x^{i-1},x_{0,x}^{i-1},y_{0,x}^{i-1},m\right)$, and set $(x_{0,x}^{i},y_{0,x}^{i})=(x_i^c,y_{x,i}^c)$.\\
$(x_{y,i}^c,y_i^c)= \text{RobustMOGDA}\left(\delta_y^{i-1},\Phi_y^{i-1},x_{0,y}^{i-1},y_{0,y}^{i-1},m\right)$, and set $(x_{0,y}^{i},y_{0,y}^{i})=(x_{y,i}^c,y_i^c)$.
}
Set $\delta_x^T=\frac{2\delta}{27(L_{xy}^2/\mu_y+L_x+\lambda_x^{T})} \quad \text{and} \quad
\delta_y^T=\frac{2\delta}{27(L_{xy}^2/\mu_x+L_y+\lambda_y^{T})}$.\\
\Return{$
x_{T+1}^c=\text{RobustMOGDA}\left(\delta_x^{T},\Phi_{x}^{T},x_{0,x}^T,y_{0,x}^T,m\right)$\\
$y_{T+1}^c=\text{RobustMOGDA}\left(\delta_y^{T},\Phi_{y}^{T},x_{0,y}^T,y_{0,y}^T,m\right)
$}
\end{algorithm}

\begin{corollary}[Efficiency of BoostMOGDA]
\label{coro: efficiency BoostMOGDA}
Fix a target accuracy $\epsilon > 0$, and a probability of failure $p \in (0,1)$. Define the algorithm parameters:
\begin{equation*}
\begin{aligned}
    &T=\left \lceil \log_{\nu}\left( \max\left(\frac{L_{xy}^2/\mu_y+L_x}{\mu_x},\frac{L_{xy}^2/\mu_x+L_y}{\mu_y}\right)
    \right) \right \rceil,\\
    &\delta=\frac{\epsilon}{4+4T}, \quad \lambda_x^i=\mu_x \nu^i \quad \text{and} \quad \lambda_y^i=\mu_y \nu^i \quad \text{with} \quad \nu=2 \quad \forall i \in [0,T].
\end{aligned}
\end{equation*}
Then the point $(x_{T+1}^c,y_{T+1}^c)=\textnormal{BoostMOGDA}(\delta,p,T)$ satisfies
\begin{equation*}
    \mathbb{P}\left(\Delta_{\Phi}(x_{T+1}^c,y_{T+1}^c) \leq \epsilon
    \right) \geq 1-p.
\end{equation*}
\end{corollary}
Note that to evaluate the sample complexity for BoostMOGDA, we need to monitor the initial squared distance for the sequence of subproblems. Define $\Delta_{\text{in}}:=\|x_0-x^*\|^2+\|y_0-y^*\|^2$, which is also the initial squared distance for the two subproblems in iteration 0. Then, We derive the upper bound for subsequent subproblems given by the following lemma.

\begin{lemma}
\label{lemma: initial gap bound}
Denote $\Delta_{\mathrm{in},x}^i$ and $\Delta_{\mathrm{in},y}^i$ the initial squared distances for the BoostMOGDA subproblems of $\Phi_{x}^i$ and $\Phi_{y}^i$, respectively. Then if all the events $E_x^i:=\left\{\|x_i^c-x_i^*\|^2+\|y_{x,i}^c-y_{x,i}^*\|^2\leq\frac{2\delta}{\mu_x+\lambda_x^{i-1}}\right\}$ and $E_y^i:=\left\{\|x_{y,i}^c-x_{y,i}^*\|^2+\|y_{i}^c-y_{i}^*\|^2\leq\frac{2\delta}{\mu_y+\lambda_y^{i-1}}\right\}$ occur for $i=0,\dots,T$, the initial squared distance can be upper bounded by
$\Delta_{\mathrm{in}}^i:= \max \left\{\Delta_{\mathrm{in},x}^i,\Delta_{\mathrm{in},y}^i\right\} \leq \mathcal{O}\left(\frac{\kappa^4\delta}{\mu}\right),$ for $i=0,1,\cdots,T.$
\end{lemma}
\proof{Proof}
We first investigate $\Delta_{\mathrm{in},x}^i:=\|x_{i}^c-x_{i+1}^*\|^2+\|y_{x,i}^c-y_{x,i+1}^*\|^2$. Applying the two-sided bound \eqref{equ: unconstrained two-side} and Proposition \ref{prop: decomposition of f and g}, we deduce
\begin{equation*}
\begin{aligned}
    \|x_{i}^c-x_{i+1}^*\|^2
    &\leq \frac{2}{\mu_x+\lambda_x^i}\left(f^i(x_{i}^c)-f^i(x_{i+1}^*)\right)\\
    &\leq \frac{2}{\mu_x+\lambda_x^i}\left(f(x_{i}^c)-f(x^*)\right)\\
    &\leq \frac{2}{\mu_x+\lambda_x^i}\left(f^{i-1}(x_{i}^c)-f^{i-1}(x_{i}^*)+\sum_{t=0}^{i-1}\frac{\lambda_x^{t}}{2}\|x_{t}^c-x_{t}^*\|^2\right)\\
    &\leq \frac{2}{\mu_x+\lambda_x^i}\left(\frac{L_f+\lambda_x^{i-1}}{2}\|x_i^c-x_i^*\|^2+\sum_{t=0}^{i-1}\frac{\lambda_x^{t}}{2}\|x_{t}^c-x_{t}^*\|^2\right)\\
    &\leq \frac{2\delta}{\mu_x+\lambda_x^i}\left(\frac{L_f+\lambda_x^{i-1}}{\mu_x+\lambda_x^{i-1}}+\sum_{t=0}^{i-1}\frac{\lambda_x^{t}}{\mu_x+\lambda_x^{t-1}}\right)\\
    &\leq \frac{2\delta}{\mu}\left(\kappa^2+2i\right)\\
    &\leq \mathcal{O}\left(\frac{\kappa^2\delta}{\mu}\right).
\end{aligned}
\end{equation*}
Similarly, we can derive the bound for the other half as
\begin{equation*}
\begin{aligned}
    \|y_{x,i}^c-y_{x,i+1}^*\|
    &\leq \|y_{x,i}^c-y_{x,i}^*\|+\|y_{x,i}^*-y_{x,i+1}^*\|\\
    &\leq \sqrt{\frac{2\delta}{\mu_x+\lambda_{x}^{i-1}}}+\frac{L_{xy}}{\mu_y}\|x_i^*-x_{i+1}^*\|\\
    &\leq \sqrt{\frac{2\delta}{\mu_x+\lambda_{x}^{i-1}}}+\frac{L_{xy}}{\mu_y}\left(\|x_i^c-x_{i}^*\|+\|x_i^c-x_{i+1}^*\|\right)\\
    &\leq \left(1+\kappa\right)\sqrt{\frac{2\delta}{\mu_x+\lambda_{x}^{i-1}}}+\kappa\|x_i^c-x_{i+1}^*\|.
\end{aligned}
\end{equation*}
Leveraging the above upper bound for $\|x_{i}^c-x_{i+1}^*\|^2$, we conclude
\begin{equation*}
    \Delta_{\mathrm{in},x}^i:=\|x_{i}^c-x_{i+1}^*\|^2+\|y_{x,i}^c-y_{x,i+1}^*\|^2\leq \mathcal{O}\left(\frac{\kappa^4\delta}{\mu}\right).
\end{equation*}
The way to bound $\Delta_{\mathrm{in},y}^i$ is identical, and is hence omitted.
\Halmos
\endproof

Now we are ready to establish the sample complexity of BoostMOGDA. For iteration $i=0,\cdots,T$, one needs to call  $\text{MOGDA}\!\big(\frac{2\delta}{27(\mu_x+\lambda_x^{i-1})},\Phi_x^{i-1}\!\!,x_{0,x}^{i-1},y_{0,x}^{i-1}\big)$ and $\text{MOGDA}\!\big(\frac{2\delta}{27(\mu_y+\lambda_y^{i-1})},\Phi_y^{i-1}\!\!,x_{0,y}^{i-1},y_{0,y}^{i-1}\big)$ for $\big\lceil18\ln\big(\frac{2T+4}{p}\big)\big\rceil$ times, respectively. 
The corresponding sample complexities are then upper bounded by $\mathcal{O}\big(\big(\frac{L+\lambda_x^{i-1}}{\mu}\big)\ln\big(\frac{(\mu_x+\lambda_x^{i-1})\Delta_{\mathrm{in},x}^{i-1}}{\delta}\big)+\frac{\sigma^2(\mu_x+\lambda_x^{i-1})}{\mu^2\delta}\big)$ 
and $\mathcal{O}\big(\big(\frac{L+\lambda_y^{i-1}}{\mu}\big)\ln\big(\frac{(\mu_y+\lambda_y^{i-1})\Delta_{\mathrm{in},y}^{i-1}}{\delta}\big)+\frac{\sigma^2(\mu_y+\lambda_y^{i-1})}{\mu^2\delta}\big)$. In iteration $T+1$, we call $\text{MOGDA}\big(\frac{2\delta}{27(L_f+\lambda_x^{T})},\Phi_x^{T},x_{0,x}^{T},y_{0,x}^{T}\big)$ and $\text{MOGDA}\big(\frac{2\delta}{27(L_g+\lambda_y^{T})},\Phi_y^{T},x_{0,y}^{T},y_{0,y}^{T}\big)$ with $\mathcal{O}\big(\big(\kappa+\frac{\lambda_x^{T}}{\mu}\big)\ln\big(\frac{(L_f+\lambda_x^{T})\Delta_{\mathrm{in},x}^{T}}{\delta}\big)+\frac{\sigma^2(L_f+\lambda_x^{T})}{\mu^2\delta}\big)$ and $\mathcal{O}\big(\big(\kappa+\frac{\lambda_y^{T}}{\mu}\big)\ln\big(\frac{(L_g+\lambda_y^{T})\Delta_{\mathrm{in},y}^{T}}{\delta}\big)+\frac{\sigma^2(L_g+\lambda_y^{T})}{\mu^2\delta}\big)$ samples, respectively. 
According to Corollary \ref{coro: efficiency BoostMOGDA}, we deduce that the total sample complexity of BoostMOGDA is on the order of
\begin{equation*}
\mathcal{O}\left(\ln(\kappa)\ln\left(\frac{\ln(\kappa)}{p}\right)
\left(\kappa^2\ln\left(\kappa\vee\frac{\Delta_{\mathrm{in}}\ln(\kappa)}{\epsilon}\right)+\frac{\sigma^2\kappa^2\ln(\kappa)}{\mu\epsilon}\right)
    \right),
\end{equation*}
which is $\tilde{\mathcal{O}}\left(\frac{\sigma^2\kappa^2}{\mu\epsilon}\right)$ after suppressing logarithm factors. This result is comparable with Theorem 2.1. of \cite{gorbunov2022clipped}, which derives a sample complexity of $\tilde{\mathcal{O}}\left(\max\left\{\kappa,\frac{\sigma^2}{\mu^2\epsilon}\right\}\right)$ to guarantee an $\epsilon$ small squared distance with high probability. By applying the upper bound of \eqref{equ: unconstrained two-side}, they can provide a guarantee for the duality gap, while the sample complexity becomes $\tilde{\mathcal{O}}\left(\frac{\sigma^2\kappa^2}{\mu\epsilon}\right)$ and is identical to our result.

\subsection{Constrained SSP problems}
Unlike the unconstrained case, the existence of constraints only allows the direct transformation of an $\mathcal{O}(\epsilon)$ expected squared distance bound to an $\mathcal{O}(\sqrt{\epsilon})$ expected duality gap. As the existing literature for constrained SC-SC SSP problems typically provides an $\mathcal{O}(\epsilon)$ expected squared distance bound with $\tilde{\mathcal{O}}(\epsilon^{-1})$ complexity dependence on $\epsilon$. Consequently, although one may still use PB-SSP to boost the confidence of general stochastic first-order oracles, directly incorporating the stochastic first-order methods in our framework may result in suboptimal sample complexities unless one can improve the analysis of existing results and derive $\tilde{\mathcal{O}}(\epsilon^{-1})$ sample complexity bounds for obtaining an $\mathcal{O}(\epsilon)$ expected duality gap. Therefore, to obtain a tighter complexity bound, we propose to adopt the strategy of \cite{lei2021stability} and use a hybrid oracle that combines SAA and stochastic first-order method for constrained SC-SC SSP. In detail, we propose to use SAA principle to construct the empirical SSP that naturally possesses a finite-sum structure, and then use stochastic variance reduced method to approximately solve this finite-sum SSP problem. In this section, we will select the loopless stochastic variance reduced extragradient (LSVRE) method given by \cite{luo2021near} as an example.   
\begin{lemma} \textbf{\emph{(\cite[Theorem 2]{luo2021near})}}
\label{lemma: SA-oracle constrained}
Suppose Assumption \ref{assump:SC-SC} and \ref{assump:smoothness} hold. Denote $(x^*_n,y^*_n)$ the saddle point of the following (deterministic) empirical SSP constructed with $n$ samples:
\begin{equation*}
    \min_{x\in \mathcal{X}}\max_{y \in \mathcal{Y}} \hat{\Phi}_n(x,y):=\frac{1}{n} \sum_{i=1}^n\Phi_{\xi_i}(x,y).
\end{equation*}
Then LSVRE can return a solution $\left(\hat{x},\hat{y}\right)$ such that  
$\mathbb{E}[\|\hat{x}-x^*_n\|^2+\|\hat{y}-y^*_n\|^2] \leq \epsilon$
within $\mathcal{O}\left((n+\sqrt{n}\kappa)\log\left(\frac{1}{\epsilon}\right)\right)$ stochastic gradient calls.
\end{lemma}
\begin{remark}
    To accommodate the above LSVRE method, we actually need to slightly strengthen the Assumption \ref{assump:smoothness} such that the Lipschitz smoothness holds almost surely for $\nabla \Phi_\xi$.
\end{remark}
Recall Lemma \ref{lemma: SAA bounded oracle} that $\mathbb{E}\left[\Delta_{\Phi}^w(x^*_n,y^*_n)\right]$ is on the order of $\mathcal{O}\left(\frac{\ell^2}{n\mu}\right)$. Then we can derive the following result by aggregating the errors of SAA and LSVRE.
\begin{theorem}[SAA-LSVRE Oracle]
Suppose Assumption \ref{assump:SC-SC}, \ref{assump:smoothness} and \ref{assump: function Lipschitz} hold.
Denote $\left(\hat{x},\hat{y}\right)$ the solution returned by applying LSVRE method to the empirical problem $\hat{\Phi}_n$ constructed by SAA. If the sample size $n$ of SAA is on the order of $\mathcal{O}\left(\frac{\ell^2}{\epsilon\mu}\right)$ and the number of calls to the stochastic gradient for LSVRE is on the order of $\mathcal{O}\left(\left(\frac{\ell^2}{\epsilon\mu}+\frac{\kappa\ell}{\sqrt{\epsilon\mu}}\right)\log\left(\frac{\ell}{\epsilon}\right)\right)$, we have $\mathbb{E}\left[\Delta_{\Phi}^w\left(\hat{x},\hat{y}\right)\right]\leq\epsilon$ where the expectation is taken with respect to both the randomness of $\hat{\Phi}_n$ and the LSVRE method.
\end{theorem}
\proof{Proof}
Denote $(x_n^*,y_n^*)$ be the saddle point of the empirical problem $\hat{\Phi}_n$ According to the definition of $\Delta_{\Phi}^w\left(\hat{x},\hat{y}\right)$, we have
\begin{equation*}
\begin{aligned}
\mathbb{E}\left[\Delta_{\Phi}^w\left(\hat{x},\hat{y}\right)\right]
&=\mathbb{E}\left[\Phi(\hat{x},y^*)-\Phi(x^*,\hat{y})\right]\\
&=\mathbb{E}\left[\Phi(\hat{x},y^*)-\Phi(x_n^*,y^*)+\Phi(x_n^*,y^*)-\Phi(x^*,y_n^*)+\Phi(x^*,y_n^*)-\Phi(x^*,\hat{y})\right]\\
&\leq \mathbb{E}\left[\Delta_{\Phi}^w\left(x_n^*,y_n^*\right)+|\Phi(\hat{x},y^*)-\Phi(x_n^*,y^*)|+|\Phi(x^*,y_n^*)-\Phi(x^*,\hat{y})|\right]\\
&\leq \mathbb{E}\left[\Delta_{\Phi}^w\left(x_n^*,y_n^*\right)+\ell\left(\|\hat{x}-x_n^*\|+\|\hat{y}-y_n^*\|\right)\right]\\
&\leq \mathbb{E}\left[\Delta_{\Phi}^w\left(x_n^*,y_n^*\right)+\ell\sqrt{2\left(\|\hat{x}-x_n^*\|^2+\|\hat{y}-y_n^*\|^2\right)}\right]\\
&\leq \mathbb{E}\left[\Delta_{\Phi}^w\left(x_n^*,y_n^*\right)\right]+\ell\sqrt{2\mathbb{E}\left[\|\hat{x}-x_n^*\|^2+\|\hat{y}-y_n^*\|^2\right]}.
\end{aligned}
\end{equation*}
Therefore, we can set $n=\mathcal{O}\left(\frac{\ell^2}{\epsilon\mu}\right)$ such that $\mathbb{E}\left[\Delta_{\Phi}^w\left(x_n^*,y_n^*\right)\right]\leq \epsilon/2$. Further applying Lemma \ref{lemma: SA-oracle constrained} and letting $\ell\sqrt{2\mathbb{E}\left[\|\hat{x}-x_n^*\|^2+\|\hat{y}-y_n^*\|^2\right]}\leq \epsilon/2$, we reach the desirable order for the number of stochastic gradient calls of the LSVRE method.
\Halmos
\endproof
Following the same spirit of the discussion in Section \ref{sec: constrained}, it is not hard to show that we can equip the above SAA-LSVRE oracle with high-probability guarantees for constrained problems in the sense of \eqref{equ: high prob guarantee}. The SAA sample size is identical to before, while the number of calls to the stochastic gradient for LSVRE is on the order of $\tilde{\mathcal{O}}\left(\frac{\ell^2}{\epsilon\mu}+\frac{\kappa^2\ell}{\sqrt{\epsilon\mu}}\right)$.

\section{Experiments}
In this section, to validate the ability of our PB-SSP framework to boost the confidence of general SSP oracles, we present the application of our method through numerical experiments on the MDP and the stochastic matrix game examples, respectively.

\subsection{Application on MDP}
First, we present the experiment for the Markov Decision Process (MDP) problem, which is a widely discussed subject within the machine learning community. When there exists a simulator that can mimic state transition and reward generation, the task of identifying the optimal policy can be reformulated as solving a constrained C-C SSP problem.

\subsubsection{SSP formulation of MDP}
An infinite-horizon average-reward MDP is specified by a tuple $\mathcal{M}=(\mathcal{S},\mathcal{A},\mathcal{P},r)$, where $\mathcal{S}$ denotes a finite state space, and $\mathcal{A}$ denotes a finite action space. $\mathcal{P}=\{P_a\}_{a\in \mathcal{A}}$ and $r=\{r_{sa}\}_{s \in \mathcal{S},a\in \mathcal{A}}$ are the state transition matrices and reward function. Specifically, when action $a$ is employed at state $s$, the system transitions to state $s'$ with a probability of $P_a(s,s')$, yielding a random reward $\hat{r}_{sa}>0$ with an expected value of $r_{sa}$. We call $\pi: \mathcal{S} \rightarrow \mathcal{P}_{\mathcal{A}}$ a stationary policy, which maps a state $s$ to a probability distribution over $\mathcal{A}$. The ultimate goal for MDP is to find an optimal policy $\pi^*$ such that the long-term average reward is maximized, irrespective of the initial states, i.e., 
\begin{equation*}
    v^*:=\max_{\pi} \,\lim_{T \rightarrow \infty} \mathbb{E}\left[\frac{1}{T}\sum_{t=0}^{T} \hat{r}_{s_t a_t} \Big| a_t \sim \pi(s_t), s_0 = s\right],
\end{equation*}
where the expectation is taken over the trajectories, which depends on both the policy and the MDP transition dynamics. The optimal Bellman equation of the MDP can be recast into an equivalent linear programming formulation (refer to \cite{chen2018scalable}), given by:
\begin{equation}
\begin{aligned}
    &\min_{v\in \mathbb{R},x\in \mathbb{R}^{|S|}} \quad v\\
    &\text{s.t.} \quad (P_a-I)x + r_a \leq v \cdot \boldsymbol{1}, \quad a \in \mathcal{A}.\\
\end{aligned}
\label{equ: LP primal}
\end{equation}
The corresponding dual formulation is: 
\begin{equation}
\begin{aligned}
    &\max_{y \in \mathbb{R}^{|S|\times|A|}} \quad \sum_{a \in \mathcal{A}} r_a^T y_a\\
    &\text{s.t.} \quad  \sum_{a \in \mathcal{A}}(P_a^T-I)y_a = 0,\|y\|_{1,1}=1, y \geq 0.\\
\end{aligned}
\label{equ: LP dual}
\end{equation}
Here, $x$ is known as the difference-value-vector and has multiple solutions obtained by adding an arbitrary constant shift. The solution $(v^*,y^*)$ of \eqref{equ: LP primal} and \eqref{equ: LP dual} correspond to the optimal average reward and stationary state-action distribution of the MDP. Furthermore, the optimal policy $\pi^*$ can be recovered from $y^*$, i.e.,
\begin{equation*}
    \mathbb{P}[\pi^*(s)=a] = \frac{y_{sa}^*}{\sum_{a' \in \mathcal{A}} y^*_{sa'}}, \quad \forall s \in \mathcal{S}.
\end{equation*}
The above can also be equivalently expressed using a saddle point formulation:
\begin{equation}
    \min_{x \in \mathcal{X}} \max_{y \in \mathcal{Y}} \Phi(x,y):=\sum_{a \in \mathcal{A}} y_a^T (P_a - I)x +\langle r, y \rangle,
    \label{equ: MDP SSP formulation}
\end{equation}
where $\mathcal{X}:=\left\{x \in \mathbb{R}^{|\mathcal{S}|}: \|x\|_{\infty}\leq U_x \right\}$,
$\mathcal{Y}:=\left\{y \in \mathbb{R}^{|\mathcal{S}| \times |\mathcal{A}|}: y \geq 0, \|y\|_{1,1}=1 \right\}$, and $U_x$ represents a bound that can be estimated. This saddle point formulation has an advantage over the linear programming formulation due to its much simpler constraints.

In practice, neither $\mathcal{P}$ nor $r$ is known, but they can be estimated through sampling from the MDP under the generator setting, see e.g. \cite{wang2017primal}. Specifically, we generate one sample by producing a transition for every $(s, a)$ pair, i.e., $\xi:=\left\{ (s,a,s',\hat{r}_{sa}): \forall s \in \mathcal{S}, a \in \mathcal{A}, s'\sim P_a(s,\cdot) \right\}$ and defining $\mathcal{P}_{\xi}$ with $P_{\xi,a}(s,s')=1$ if $(s,a) \rightarrow s'$ is observed and $P_{\xi,a}(s,s')=0$ otherwise. This leads to
\begin{equation*}
    \Phi_{\xi}(x,y):=\sum_{a \in \mathcal{A}} y_a^T (P_{\xi,a} - I)x +\langle \hat{r}, y \rangle,
    \label{equ: MDP SSP sample}
\end{equation*}
by direct computation, we know $\Phi(x,y)=\mathbb{E}_{\xi}[\Phi_{\xi}(x,y)]$. Under this setting, a standard SAA approach can be adopted as the in-expectation oracle that returns a policy with $\epsilon$-optimal average reward in expectation, as discussed in \cite{zhang2021generalization}. 

\subsubsection{Numerical experiments}
As discussed in section \ref{sec: constrained C-C}, it is essential to incorporate quadratic regularization terms into the constrained C-C SSP problem to make our PB-SSP applicable. Consequently, we formulate the following regularized objective function:
\begin{equation*}
\begin{aligned}
    \Phi_{\alpha}(x,y):=&\Phi(x,y)+\frac{\epsilon}{4D_x^2}\|x\|^2-\frac{\epsilon}{4D_y^2}\left\|y-\boldsymbol{1}/(|S|\cdot|A|)\right\|^2.
\end{aligned}
\end{equation*}
In our experiment, we set $|\mathcal{S}|=100$, $|\mathcal{A}|=10$, and randomly generate $0<r_{sa}<1$ and $\{P_a\}_{a \in \mathcal{A}}$. To demonstrate the efficacy of our procedure, we suppose that $\hat{r}_{sa}$ follows heavy-tailed gamma distributions with a universal variance $\text{Var}[\hat{r}_{sa}]:= \sigma^2_r = 1$. We estimate $U_x$ to be $0.5$ in our example, which gives us $D_x^2=25$. We also have $D_y^2=2$ since $\mathcal{Y}$ is the unit simplex. We aim for a target accuracy and a probability of failure of $\epsilon=0.01$ and $p=1\%$, respectively. 

For comparison purposes, we also evaluate SAA \cite{zhang2021generalization} and the RDE approach (as mentioned in section \ref{subsec: constrained SAA}). All the empirical proximal subproblems of IPPA are solved by the proximal extra gradient method. The parameter settings, informed by the preceding corollaries and discussions, are theoretically sound but tend to be excessive for practical applications. To facilitate the numerical experiment, we follow the general framework of the procedures and select the parameters independently. To ensure fairness, we first standardize the SAA with a sample size of $10^5$, serving as the basic oracle operation for all procedures, irrespective of the values of other parameters. The sample size for the gradient oracles (when required) is consistently set at $10^4$. For the RDE approach, only the parameter $m$ needs to be determined. In our PB-SSP procedure, we always set $\nu=4$ as the base number as it is consistently the best performing parameter in the experiments, the users have the flexibility to adjust  $T$ and $m$. Table \ref{table: MDP results quadratic} summarizes the performance results along with different parameter settings. We document the effort expended in invoking the basic SAA oracle and the resultant duality gap for the original SSP problem \eqref{equ: MDP SSP formulation}, with statistics estimated from $1,000$ independent macro replications. The decimal points in the column labeled ``\# of calls'' arise from the samples used by the gradient oracle, where we approximate the effort for one call of the gradient oracle as being $10\%$ that of the basic SAA oracle.

\begin{table}[htb!]
\centering
\caption{Comparison of procedures on the SSP problem for MDP (quadratic regularization term).}
\begin{tabular}{C{2.8cm} | C{1.0cm}  C{1.0cm}  C{1.0cm} | C{2.0cm}  C{2.5cm} C{3.5cm}}
\hline
     & $\nu$ & T & $m$ & $\#$ of calls & $\mathbb{E}\left[\Delta_{\Phi}(\bar{x},\bar{y})\right]$ & $\mathbb{P}\left[\Delta_{\Phi}(\bar{x},\bar{y}) > 0.01 \right]$\\
\hline
    SAA & - & - & - & 1 & 0.0093 & $26.5\%$\\
\hline
    \multirow{4}{*}{SAA+RDE} 
    & - & - & 3 & 3.6 & 0.0091 & $23.3\%$\\
    & - & - & 9 & 10.8 & 0.0086 & $13.5\%$\\
    & - & - & 99 & 118.8 & 0.0080 & $5.4\%$\\
    & - & - & 499 & 598.8 & 0.0078 & $2.6\%$\\
\hline
    \multirow{4}{*}{SAA+PB-SSP}
    & 4 & 7 & 3 & 54.6 & 0.0076 & $2.1\%$\\
    & 4 & 7 & 5 & 91.0 & 0.0074 & $1.3\%$\\
    & 4 & 8 & 3 & \boxit{6.6 cm}60.6 & 0.0069 & $0.2\%$\\    
    & 4 & 8 & 5 & 101.0 & 0.0066 & $0\%$\\
\hline 
\end{tabular}  
\label{table: MDP results quadratic}
\end{table}

In addition, since $\mathcal{Y}$ is the unit simplex in this example, it allows us to integrate an entropy regularization term for $y$, which exhibits strong concavity in the $\ell^1$ norm. Both theoretical and empirical studies have shown entropy regularization terms to deliver enhanced performance for optimization problems in a simplex. As such, we devise the following alternative regularized objective function:
\begin{equation*}
\begin{aligned}
    \Phi_{\alpha}(x,y):=&\Phi(x,y)+\frac{\epsilon}{4D_x^2}\|x\|^2-\frac{\epsilon}{4\log(|\mathcal{S}||\mathcal{A}|)}\sum_{s,a}y_{sa}\log(y_{sa}).
\end{aligned}
\end{equation*}
Likewise, the proximal point term in $\Phi_y^i$ is also substituted with the Kullback–Leibler divergence $D_{\text{KL}}(y|y_i^c)$ accordingly. It is noteworthy that our analysis for PB-SSP could actually be extended to SC-SC SSP problems in general norms. Table \ref{table: MDP results entropy} presents the results obtained after the integration of the entropy regularization term.
\begin{table}[htb!]
\centering
\caption{Comparison of procedures on the SSP problem for MDP (entropy regularization term).}
\begin{tabular}{C{2.8cm} | C{1.0cm}  C{1.0cm}  C{1.0cm} | C{2.0cm}  C{2.5cm} C{3.5cm}}
\hline
     & $\nu$ & T & $m$ & $\#$ of calls & $\mathbb{E}\left[\Delta_{\Phi}(\bar{x},\bar{y})\right]$ & $\mathbb{P}\left[\Delta_{\Phi}(\bar{x},\bar{y}) > 0.01 \right]$\\
\hline
    SAA & - & - & - & 1 & 0.0092 & $25.0\%$\\
\hline
    \multirow{4}{*}{SAA+RDE} 
    & - & - & 3 & 3.6 & 0.0090 & $21.2\%$\\
    & - & - & 9 & 10.8 & 0.0086 & $13.1\%$\\
    & - & - & 99 & 118.8 & 0.0080 & $4.6\%$\\
    & - & - & 499 & 598.8 & 0.0077 & $3.1\%$\\
\hline
    \multirow{4}{*}{SAA+PB-SSP}
    & 4 & 7 & 3 & \boxit{6.6 cm}54.6 & 0.0073 & $1.0\%$\\
    & 4 & 7 & 5 & 91.0 & 0.0071 & $0.9\%$\\
    & 4 & 8 & 3 & 60.6 & 0.0064 & $0.2\%$\\
    & 4 & 8 & 5 & 101.0 & 0.0061 & $0\%$\\
\hline 
\end{tabular}  
\label{table: MDP results entropy}
\end{table}

Upon examining the two tables, we observe that adding the entropy regularization term generally yields relatively superior results compared to adding the quadratic one. An intuitive explanation is that the entropy regularization term more accurately reflects the geometric properties of the simplex domain. Furthermore, the three procedures exhibit similar performance in both scenarios. The basic SAA oracle can produce solutions with small duality gaps in expectation but lacks high-probability guarantees. The RDE approach can incrementally enhance confidence as $m$ goes up, yet its marginal benefits diminish rapidly. Even after 598.8 calls to the basic SAA oracle, it falls short of achieving our target probability of failure.
In contrast, the PB-SSP framework demonstrates a significant advantage by additionally integrating the inexact proximal point method. In general, for a fixed base number $\nu$, larger values of $T$ and $m$ yield improved results, and a judicious selection of these three parameters can achieve our target high-probability guarantees at a reasonable cost.
When the quadratic regularization term is added, PB-SSP requires 60.6 calls to the basic SAA oracle, whereas only 54.6 calls are needed when the entropy regularization term is employed. 

\subsection{Application on stochastic matrix game}
In this subsection, we consider the two-player stochastic matrix game given by
\begin{equation*}
\min_{x\in\mathcal{X}}\max_{y\in\mathcal{Y}} \Phi(x,y):=x^T\mathbb{E}\left[A_{\xi}\right]y,
\end{equation*}
where $\mathcal{X}:=\{x\in\mathbb{R}^{N_x}:x\geq0,\mathbf{1}^Tx=1\}$ and $\mathcal{Y}:=\{y\in\mathbb{R}^{N_y}:y\geq0,\mathbf{1}^Ty=1\}$. The decision variables $x$ and $y$ are the mixed strategies of the two players, and $A_{\xi}\in \mathbb{R}^{N_x\times N_y}$ is the stochastic payoff. To apply the PB-SSP framework, we consider the following entropy-regularized problem:
\begin{equation*}
\min_{x\in\mathcal{X}}\max_{y\in\mathcal{Y}} \Phi_{\alpha}(x,y):=x^T\mathbb{E}\left[A_{\xi}\right]y + \frac{\epsilon}{4\log(N_x)}\sum_{i=1}^{N_x}x_i\log(x_i)-\frac{\epsilon}{4\log(N_y)}\sum_{i=1}^{N_y}y_i\log(y_i),
\end{equation*}
since entropy regularization terms are reported to deliver better performance in the MDP example.
In our experimental setup, we specify the dimensions of the strategy spaces with $N_x=100$ for the first player and $N_y=200$ for the second player. The elements $a_{ij}$ of the expected payoff matrix $A:=\mathbb{E}[A_{\xi}]$ are generated randomly within the range of 0 to 1. Likewise, we model the $\hat{a}_{ij}$, representing the stochastic components of the payoff, to follow Gamma distributions with a universal variance of $\sigma^2_A=1$, making sub-Gaussian assumption not applicable. The target accuracy and probability of failure are still $\epsilon=0.01$ and $p=1\%$. Next, we construct two basic SSP oracles for comparison. The first one is by constructing empirical problems by SAA with a sample size of $5,000$ and then solving by the proximal extra gradient method. The second oracle is the stochastic proximal extra gradient (SPEG) algorithm whose iterates are given by,
\begin{equation*}
\begin{aligned}
&\tilde{x}^{t}=\argmin_{x\in\mathcal{X}} \Phi_{\bar{\xi}^t_1}(x,y^{t-1})+\eta_{t}\sum_{i=1}^{N_x}x_i\log\left(\frac{x_i}{x_i^{t-1}}\right),
&&\tilde{y}^{t}=\argmax_{y\in\mathcal{Y}} \Phi_{\bar{\xi}^t_1}(x^{t-1},y)+\eta_{t}\sum_{i=1}^{N_y}y_i\log\left(\frac{y_i}{y_i^{t-1}}\right),\\
&x^{t}=\argmin_{x\in\mathcal{X}} \Phi_{\bar{\xi}^t_2}(x,y^{t-1})+\eta_{t}\sum_{i=1}^{N_x}x_i\log\left(\frac{x_i}{\tilde{x}_i^{t}}\right),
&&y^{t}=\argmax_{y\in\mathcal{Y}} \Phi_{\bar{\xi}^t_2}(x^{t-1},y)+\eta_{t}\sum_{i=1}^{N_y}y_i\log\left(\frac{y_i}{\tilde{y}_i^{t}}\right).
\end{aligned}
\end{equation*}
As the cross term of this problem is bilinear, we are able to maintain the linearization of the cross as itself in the proximal gradient step, which makes the proximal gradient step identical to a proximal point step. 
We fix the iteration count at 2,000, calibrate an appropriate constant stepsize $\eta_t$, and adopt independent samples $(\bar{\xi}_1^t,\bar{\xi}_2^t)$ with a batch size of 10. The average of all iterates is returned as the solution.  
The numerical results are documented in Table \ref{table: Matrix results entropy} and \ref{table: Matrix results entropy first order}.

\begin{table}[htb!]
\centering
\caption{Comparison of SAA procedures on the SSP problem for stochastic matrix game.}
\begin{tabular}{C{2.8cm} | C{1.0cm}  C{1.0cm}  C{1.0cm} | C{2.0cm}  C{2.5cm} C{3.5cm}}
\hline
     & $\nu$ & T & $m$ & $\#$ of calls & $\mathbb{E}\left[\Delta_{\Phi}(\bar{x},\bar{y})\right]$ & $\mathbb{P}\left[\Delta_{\Phi}(\bar{x},\bar{y}) > 0.01 \right]$\\
\hline
    \multirow{1}{*}{SAA} 
    & - & - & - & 1 & 0.0100 & $45.9\%$\\
\hline
    \multirow{4}{*}{SAA+RDE} 
    & - & - & 3 & 3.6 & 0.0099 & $43.7\%$\\
    & - & - & 9 & 10.8 & 0.0096 & $35.5\%$\\
    & - & - & 99 & 118.8 & 0.0093 & $27.3\%$\\
    & - & - & 499 & 598.8 & 0.0091 & $24.1\%$\\
\hline
    \multirow{4}{*}{SAA+PB-SSP}
    & 4 & 5 & 3 & 42.6 & 0.0076 & $2.8\%$\\
    & 4 & 5 & 5 & 71.0 & 0.0074 & $2.3\%$\\
    & 4 & 6 & 3 & 48.6 & 0.0068 & $2.0\%$\\
    & 4 & 6 & 5 & \boxit{6.6 cm}81.0 & 0.0064 & $0\%$\\
\hline
\end{tabular}  
\label{table: Matrix results entropy}
\end{table}
\begin{table}[htb!]
\centering
\caption{Comparison of first-order procedures on the SSP problem for stochastic matrix game.}
\begin{tabular}{C{2.8cm} | C{1.0cm}  C{1.0cm}  C{1.0cm} | C{2.0cm}  C{2.5cm} C{3.5cm}}
\hline
     & $\nu$ & T & $m$ & $\#$ of calls & $\mathbb{E}\left[\Delta_{\Phi}(\bar{x},\bar{y})\right]$ & $\mathbb{P}\left[\Delta_{\Phi}(\bar{x},\bar{y}) > 0.01 \right]$\\
\hline
    \multirow{1}{*}{SPEG} 
    & - & - & - & 1 & 0.0095& $27.8\%$\\
\hline
    \multirow{4}{*}{SPEG+RDE} 
    & - & - & 3 & 3.6 & 0.0095 & $26.3\%$\\
    & - & - & 9 & 10.8 & 0.0093 & $21.3\%$\\
    & - & - & 99 & 118.8 & 0.0090 & $12.5\%$\\
    & - & - & 499 & 598.8 & 0.0089 & $11.8\%$\\
\hline
    \multirow{4}{*}{SPEG+PB-SSP}
    & 4 & 1 & 3 & 18.6 & 0.0085 & $4.5\%$\\
    & 4 & 1 & 5 & 31.0 & 0.0085 & $4.9\%$\\
    & 4 & 2 & 3 & \boxit{6.6 cm}24.6 & 0.0075 & $0.3\%$\\
    & 4 & 2 & 5 & 41.0 & 0.0075 & $0.4\%$\\
\hline
\end{tabular}  
\label{table: Matrix results entropy first order}
\end{table}

Similar to the MDP example, the RDE approach alone is hard to boost SSP oracles into high confidence, no matter for SAA or first-order procedures. By contrast, the PB-SSP framework can obtain a high-confidence solution within tens of invocations of the given basic oracle.




\bibliographystyle{plain} 
\bibliography{main} 

\begin{thebibliography}{10}

\bibitem{alacaoglu2022stochastic}
Ahmet Alacaoglu and Yura Malitsky.
\newblock Stochastic variance reduction for variational inequality methods.
\newblock In {\em Conference on Learning Theory}, pages 778--816. PMLR, 2022.

\bibitem{bach2019universal}
Francis Bach and Kfir~Y Levy.
\newblock A universal algorithm for variational inequalities adaptive to smoothness and noise.
\newblock In {\em Conference on learning theory}, pages 164--194. PMLR, 2019.

\bibitem{beznosikov2020gradient}
Aleksandr Beznosikov, Abdurakhmon Sadiev, and Alexander Gasnikov.
\newblock Gradient-free methods with inexact oracle for convex-concave stochastic saddle-point problem.
\newblock In {\em International Conference on Mathematical Optimization Theory and Operations Research}, pages 105--119. Springer, 2020.

\bibitem{chavdarova2019reducing}
Tatjana Chavdarova, Gauthier Gidel, Fran{\c{c}}ois Fleuret, and Simon Lacoste-Julien.
\newblock Reducing noise in gan training with variance reduced extragradient.
\newblock {\em Advances in Neural Information Processing Systems}, 32, 2019.

\bibitem{chen2018scalable}
Yichen Chen, Lihong Li, and Mengdi Wang.
\newblock Scalable bilinear pi learning using state and action features.
\newblock In {\em International Conference on Machine Learning}, pages 834--843. PMLR, 2018.

\bibitem{chen2014optimal}
Yunmei Chen, Guanghui Lan, and Yuyuan Ouyang.
\newblock Optimal primal-dual methods for a class of saddle point problems.
\newblock {\em SIAM Journal on Optimization}, 24(4):1779--1814, 2014.

\bibitem{chen2017accelerated}
Yunmei Chen, Guanghui Lan, and Yuyuan Ouyang.
\newblock Accelerated schemes for a class of variational inequalities.
\newblock {\em Mathematical Programming}, 165(1):113--149, 2017.

\bibitem{davis2021low}
Damek Davis, Dmitriy Drusvyatskiy, Lin Xiao, and Junyu Zhang.
\newblock From low probability to high confidence in stochastic convex optimization.
\newblock {\em Journal of Machine Learning Research}, 22(49):1--38, 2021.

\bibitem{fallah2020optimal}
Alireza Fallah, Asuman Ozdaglar, and Sarath Pattathil.
\newblock An optimal multistage stochastic gradient method for minimax problems.
\newblock In {\em 2020 59th IEEE Conference on Decision and Control (CDC)}, pages 3573--3579. IEEE, 2020.

\bibitem{farnia2021train}
Farzan Farnia and Asuman Ozdaglar.
\newblock Train simultaneously, generalize better: Stability of gradient-based minimax learners.
\newblock In {\em International Conference on Machine Learning}, pages 3174--3185. PMLR, 2021.

\bibitem{gidel2018variational}
Gauthier Gidel, Hugo Berard, Ga{\"e}tan Vignoud, Pascal Vincent, and Simon Lacoste-Julien.
\newblock A variational inequality perspective on generative adversarial networks.
\newblock {\em arXiv preprint arXiv:1802.10551}, 2018.

\bibitem{goodfellow2020generative}
Ian Goodfellow, Jean Pouget-Abadie, Mehdi Mirza, Bing Xu, David Warde-Farley, Sherjil Ozair, Aaron Courville, and Yoshua Bengio.
\newblock Generative adversarial networks.
\newblock {\em Communications of the ACM}, 63(11):139--144, 2020.

\bibitem{gorbunov2022clipped}
Eduard Gorbunov, Marina Danilova, David Dobre, Pavel Dvurechenskii, Alexander Gasnikov, and Gauthier Gidel.
\newblock Clipped stochastic methods for variational inequalities with heavy-tailed noise.
\newblock {\em Advances in Neural Information Processing Systems}, 35:31319--31332, 2022.

\bibitem{hsieh2019convergence}
Yu-Guan Hsieh, Franck Iutzeler, J{\'e}r{\^o}me Malick, and Panayotis Mertikopoulos.
\newblock On the convergence of single-call stochastic extra-gradient methods.
\newblock {\em Advances in Neural Information Processing Systems}, 32, 2019.

\bibitem{hsieh2020explore}
Yu-Guan Hsieh, Franck Iutzeler, J{\'e}r{\^o}me Malick, and Panayotis Mertikopoulos.
\newblock Explore aggressively, update conservatively: Stochastic extragradient methods with variable stepsize scaling.
\newblock {\em Advances in Neural Information Processing Systems}, 33:16223--16234, 2020.

\bibitem{hsu2016loss}
Daniel Hsu and Sivan Sabato.
\newblock Loss minimization and parameter estimation with heavy tails.
\newblock {\em Journal of Machine Learning Research}, 17(1):543--582, 2016.

\bibitem{huang2022new}
Kevin Huang and Shuzhong Zhang.
\newblock New first-order algorithms for stochastic variational inequalities.
\newblock {\em SIAM Journal on Optimization}, 32(4):2745--2772, 2022.

\bibitem{iusem2017extragradient}
Alfredo~N Iusem, Alejandro Jofr{\'e}, Roberto~Imbuzeiro Oliveira, and Philip Thompson.
\newblock Extragradient method with variance reduction for stochastic variational inequalities.
\newblock {\em SIAM Journal on Optimization}, 27(2):686--724, 2017.

\bibitem{juditsky2011solving}
Anatoli Juditsky, Arkadi Nemirovski, and Claire Tauvel.
\newblock Solving variational inequalities with stochastic mirror-prox algorithm.
\newblock {\em Stochastic Systems}, 1(1):17--58, 2011.

\bibitem{laguel2023high}
Yassine Laguel, Necdet~Serhat Aybat, and Mert G{\"u}rb{\"u}zbalaban.
\newblock High probability and risk-averse guarantees for a stochastic accelerated primal-dual method.
\newblock {\em arXiv preprint arXiv:2304.00444}, 2023.

\bibitem{lei2021stability}
Yunwen Lei, Zhenhuan Yang, Tianbao Yang, and Yiming Ying.
\newblock Stability and generalization of stochastic gradient methods for minimax problems.
\newblock In {\em International Conference on Machine Learning}, pages 6175--6186. PMLR, 2021.

\bibitem{luo2019stochastic}
Luo Luo, Cheng Chen, Yujun Li, Guangzeng Xie, and Zhihua Zhang.
\newblock A stochastic proximal point algorithm for saddle-point problems.
\newblock {\em arXiv preprint arXiv:1909.06946}, 2019.

\bibitem{luo2021near}
Luo Luo, Guangzeng Xie, Tong Zhang, and Zhihua Zhang.
\newblock Near optimal stochastic algorithms for finite-sum unbalanced convex-concave minimax optimization.
\newblock {\em arXiv preprint arXiv:2106.01761}, 2021.

\bibitem{mishchenko2020revisiting}
Konstantin Mishchenko, Dmitry Kovalev, Egor Shulgin, Peter Richt{\'a}rik, and Yura Malitsky.
\newblock Revisiting stochastic extragradient.
\newblock In {\em International Conference on Artificial Intelligence and Statistics}, pages 4573--4582. PMLR, 2020.

\bibitem{namkoong2017variance}
Hongseok Namkoong and John~C Duchi.
\newblock Variance-based regularization with convex objectives.
\newblock {\em Advances in neural information processing systems}, 30, 2017.

\bibitem{natole2018stochastic}
Michael Natole, Yiming Ying, and Siwei Lyu.
\newblock Stochastic proximal algorithms for auc maximization.
\newblock In {\em International Conference on Machine Learning}, pages 3710--3719. PMLR, 2018.

\bibitem{nemirovski2009robust}
Arkadi Nemirovski, Anatoli Juditsky, Guanghui Lan, and Alexander Shapiro.
\newblock Robust stochastic approximation approach to stochastic programming.
\newblock {\em SIAM Journal on optimization}, 19(4):1574--1609, 2009.

\bibitem{nemirovski2002efficient}
Arkadi Nemirovski and Reuven~Y Rubinstein.
\newblock An efficient stochastic approximation algorithm for stochastic saddle point problems.
\newblock In {\em Modeling Uncertainty}, pages 156--184. Springer, 2002.

\bibitem{nemirovskij1983problem}
Arkadij~Semenovi{\v{c}} Nemirovskij and David~Borisovich Yudin.
\newblock {\em Problem complexity and method efficiency in optimization}.
\newblock John Wiley \& Sons, Inc., New York, 1983.

\bibitem{nesterov2005smooth}
Yu~Nesterov.
\newblock Smooth minimization of non-smooth functions.
\newblock {\em Mathematical programming}, 103:127--152, 2005.

\bibitem{nouiehed2019solving}
Maher Nouiehed, Maziar Sanjabi, Tianjian Huang, Jason~D Lee, and Meisam Razaviyayn.
\newblock Solving a class of non-convex min-max games using iterative first order methods.
\newblock In H.~Wallach, H.~Larochelle, A.~Beygelzimer, F.~d\textquotesingle Alch\'{e}-Buc, E.~Fox, and R.~Garnett, editors, {\em Advances in Neural Information Processing Systems}, volume~32. Curran Associates, Inc., 2019.

\bibitem{palaniappan2016stochastic}
Balamurugan Palaniappan and Francis Bach.
\newblock Stochastic variance reduction methods for saddle-point problems.
\newblock In {\em Advances in Neural Information Processing Systems}, volume~29, pages 1416--1424, 2016.

\bibitem{polyak1990new}
Boris~T Polyak.
\newblock New stochastic approximation type procedures.
\newblock {\em Automat. i Telemekh}, 7(98-107):2, 1990.

\bibitem{puterman2014markov}
Martin~L Puterman.
\newblock {\em Markov decision processes: discrete stochastic dynamic programming}.
\newblock John Wiley \& Sons, 2014.

\bibitem{roughgarden2010algorithmic}
Tim Roughgarden.
\newblock Algorithmic game theory.
\newblock {\em Communications of the ACM}, 53(7):78--86, 2010.

\bibitem{shalev2013stochastic}
Shai Shalev-Shwartz and Tong Zhang.
\newblock Stochastic dual coordinate ascent methods for regularized loss minimization.
\newblock {\em Journal of Machine Learning Research}, 14(Feb):567--599, 2013.

\bibitem{sinha2017certifying}
Aman Sinha, Hongseok Namkoong, Riccardo Volpi, and John Duchi.
\newblock Certifying some distributional robustness with principled adversarial training.
\newblock {\em arXiv preprint arXiv:1710.10571}, 2017.

\bibitem{von2007theory}
John Von~Neumann and Oskar Morgenstern.
\newblock Theory of games and economic behavior.
\newblock In {\em Theory of games and economic behavior}. Princeton university press, 2007.

\bibitem{wang2017primal}
Mengdi Wang.
\newblock Primal-dual pi learning: Sample complexity and sublinear run time for ergodic markov decision problems.
\newblock {\em arXiv preprint arXiv:1710.06100}, 2017.

\bibitem{xu2010sample}
Huifu Xu.
\newblock Sample average approximation methods for a class of stochastic variational inequality problems.
\newblock {\em Asia-Pacific Journal of Operational Research}, 27(01):103--119, 2010.

\bibitem{yan2020optimal}
Yan Yan, Yi~Xu, Qihang Lin, Wei Liu, and Tianbao Yang.
\newblock Optimal epoch stochastic gradient descent ascent methods for min-max optimization.
\newblock In {\em Advances in Neural Information Processing Systems}, pages 5789--5800. Curran Associates, Inc., 2020.

\bibitem{yan2019stochastic}
Yan Yan, Yi~Xu, Qihang Lin, Lijun Zhang, and Tianbao Yang.
\newblock Stochastic primal-dual algorithms with faster convergence than ${O}(1/\sqrt{T})$ for problems without bilinear structure.
\newblock {\em arXiv preprint arXiv:1904.10112}, 2019.

\bibitem{zhang2021generalization}
Junyu Zhang, Mingyi Hong, Mengdi Wang, and Shuzhong Zhang.
\newblock Generalization bounds for stochastic saddle point problems.
\newblock In {\em International Conference on Artificial Intelligence and Statistics}, pages 568--576. PMLR, 2021.

\bibitem{zhang2017stochastic}
Yuchen Zhang and Lin Xiao.
\newblock Stochastic primal-dual coordinate method for regularized empirical risk minimization.
\newblock {\em The Journal of Machine Learning Research}, 18(1):2939--2980, 2017.

\bibitem{zhao2022accelerated}
Renbo Zhao.
\newblock Accelerated stochastic algorithms for convex-concave saddle-point problems.
\newblock {\em Mathematics of Operations Research}, 47(2):1443--1473, 2022.

\end{thebibliography}



%

\end{document}